\documentclass[a4,11pt]{amsart}
\usepackage{amsthm,amssymb,amsmath,color,epsfig}

\newcommand{\red}{\color{darkred}}
\newcommand{\blue}{\color{darkblue}}
\newcommand{\green}{\color{darkgreen}}

\newcommand{\clt}{central limit theorem}

\newcommand{\ex}{{\rm e}\,}

\newcommand{\asy}{asymptotic}

\newcommand{\ts}{time series}

\definecolor{darkblue}{rgb}{.1, 0.1,.8}
\definecolor{darkgreen}{rgb}{0,0.8,0.2}
\definecolor{darkred}{rgb}{.8, .1,.1}

\newtheorem{lemma}{Lemma}[section]

\newtheorem{theorem}[lemma]{Theorem}

\newtheorem{proposition}[lemma]{Proposition}
\newtheorem{definition}[lemma]{Definition}
\newtheorem{corollary}[lemma]{Corollary}
\newtheorem{example}[lemma]{Example}
\newtheorem{exercise}[lemma]{Exercise}
\newtheorem{remark}[lemma]{Remark}
\newtheorem{fig}[lemma]{Figure}
\newtheorem{tab}[lemma]{Table}

\newcommand{\MC}{Markov chain}

\newcommand{\bfu}{{\bf u}}
\newcommand{\bfv}{{\bf v}}

\newcommand{\bth}{\begin{theorem}}
\newcommand{\ethe}{\end{theorem}}

\newcommand{\bre}{\begin{remark}\em }
\newcommand{\ere}{\end{remark}}

\newcommand{\ble}{\begin{lemma}}
\newcommand{\ele}{\end{lemma}}
\newcommand{\sre}{stochastic recurrence equation}

\newcommand{\bde}{\begin{definition}}
\newcommand{\ede}{\end{definition}}
\newcommand{\bco}{\begin{corollary}}
\newcommand{\eco}{\end{corollary}}

\newcommand{\bpr}{\begin{proposition}}
\newcommand{\epr}{\end{proposition}}

\newcommand{\bexer}{\begin{exercise}}
\newcommand{\eexer}{\end{exercise}}

\newcommand{\bexam}{\begin{example}}
\newcommand{\eexam}{\end{example}}

\newcommand{\bfi}{\begin{fig}}
\newcommand{\efi}{\end{fig}}

\newcommand{\btab}{\begin{tab}}
\newcommand{\etab}{\end{tab}}

\newcommand{\rv}{random variable}

\newcommand{\var}{{\rm var}}

\newcommand{\as}{{\rm a.s.}}

\newcommand{\rhs}{right-hand side}
\newcommand{\df}{distribution function}

\newcommand{\beao}{\begin{eqnarray*}}
\newcommand{\eeao}{\end{eqnarray*}\noindent}

\newcommand{\beam}{\begin{eqnarray}}
\newcommand{\eeam}{\end{eqnarray}\noindent}

\newcommand{\beqq}{\begin{equation}}
\newcommand{\eeqq}{\end{equation}\noindent}

\newcommand{\bce}{\begin{center}}
\newcommand{\ece}{\end{center}}

\newcommand{\barr}{\begin{array}}
\newcommand{\earr}{\end{array}}

\newcommand{\std}{\stackrel{d}{\rightarrow}}
\newcommand{\stas}{\stackrel{\rm a.s.}{\rightarrow}}

\newcommand{\stv}{\stackrel{v}{\rightarrow}}

\newcommand{\eqd}{\stackrel{d}{=}}

\newcommand{\vague}{\stackrel{\lower0.2ex\hbox{$\scriptscriptstyle
                    \it{v} $}}{\rightarrow}}
\newcommand{\weak}{\stackrel{\lower0.2ex\hbox{$\scriptscriptstyle
                    \it{w} $}}{\rightarrow}}
\newcommand{\what}{\stackrel{\lower0.2ex\hbox{$\scriptscriptstyle
                    \it{\hat{w}} $}}{\rightarrow}}

\newcommand{\bdis}{\begin{displaymath}}
\newcommand{\edis}{\end{displaymath}\noindent}

\renewcommand{\P}{\mathbb P}
\newcommand{\R}{\mathbb{R}}

\newcommand{\nto}{n\to\infty}

\newcommand{\xto}{x\to\infty}

\newcommand{\ov}{\overline}
\newcommand{\wt}{\widetilde}

\newcommand{\vep}{\varepsilon}

\newcommand{\la}{\lambda}

\newcommand{\regvary}{regularly varying}
\newcommand{\slvary}{slowly varying}
\newcommand{\regvar}{regular variation}

\newcommand{\bbr}{{\mathbb R}}

\newcommand{\bbz}{{\mathbb Z}}

\newcommand{\bbf}{{\mathcal F}}
\newcommand{\bbs}{{\mathbb S}}

\newcommand{\con}{convergence}

\newcommand{\st}{such that}
\newcommand{\fif}{if and only if}
\newcommand{\wrt}{with respect to}

\newcommand{\fct}{function}

\newcommand{\ds}{distribution}

\newcommand{\rep}{representation}

\newcommand{\seq}{sequence}

\newcommand{\pro}{probabilit}

\newcommand{\ms}{measure}

\newcommand{\bfx}{{\bf x}}
\newcommand{\bfX}{{\bf X}}

\newcommand{\bfY}{{\bf Y}}
\newcommand{\bfy}{{\bf y}}

\newcommand{\bfZ}{{\bf Z}}

\newcommand{\bfs}{{\bf s}}

\def\1{\ensuremath{\mathrm{1}\hspace{-.35em} \mathrm{1}}} 

\def\E{{\mathbb E}}

\def\P{{\mathbb{P}}}
\def\R{\mathbb{R}}

\renewcommand{\le}{\ensuremath{\leqslant}}
\renewcommand{\ge}{\ensuremath{\geqslant}}

\newcommand{\introo}[2]{{\left]{#1,\,#2\,}\right[\kern1pt}}

\newcommand{\intrfo}[2]{{\left[{#1,\,#2}\right[\kern1pt}}

\begin{document}

\title{Heavy tails for an alternative stochastic perpetuity model}



\date{}
\today

\thanks{Thomas Mikosch's and Mohsen Rezapour's research is partly supported
by the Danish Research Council Grant DFF-4002-00435. Financial support by the ANR network AMERISKA ANR 14 CE20 0006 01 is gratefully
acknowledged by Thomas Mikosch and Olivier Wintenberger.}
\author[T. Mikosch]{Thomas Mikosch}
\address{Department  of Mathematics\\
University of Copenhagen\\
Universitetsparken 5\\
DK-2100 Copenhagen\\
Denmark}
\email{mikosch@math.ku.dk}
\author[M. Rezapour]{Mohsen Rezapour}
\address{Department of Statistics\\ 
Faculty of Mathematics and Computer Science\\ 
Shahid Bahonar University of Kerman\\ 
Kerman\\ Iran}
\email{mohsenrzp@gmail.com}
\author[O. Wintenberger]{Olivier Wintenberger}
\address{Sorbonne Universit\'es\\
UPMC Universit\'e Paris 06\\ 
F-75005, Paris\\ France}
\email{olivier.wintenberger@upmc.fr}

\begin{abstract}
In this paper we consider a stochastic model of perpetuity-type. In contrast to the classical
affine perpetuity model of Kesten \cite{kesten:1973} and Goldie \cite{goldie:1991} all discount factors in the model
are mutually independent. We prove that the tails of the \ds\ of this model are \regvary\ both in the 
univariate and multivariate cases. Due to the additional randomness in the model the tails are not pure power laws
as in the Kesten-Goldie setting but involve a logarithmic term. \end{abstract}


\maketitle
\section{Problem description}\setcounter{equation}{0}
Consider an array $(X_{ni})$ of iid \rv s with generic \seq\ $(X_i)=(X_{1i})$ and $X=X_1$.
We define a stochastic {\em perpetuity} in the following way:
\beao
\wt Y_1&=& X_{11}\,,\\
\wt Y_2&=& X_{11}X_{12}+X_{21}\,,\\
\wt Y_3&=& X_{11}X_{12}X_{13}+X_{21}X_{22}+X_{31}\,,\ldots\,.
\eeao
At any time $i$, each of the investments in the previous and current periods $j=1,\ldots,i$ 
gets discounted by an independent factor $X_{ij}$. 
Therefore $(\wt Y_n)$ can be interpreted as the dynamics of a perpetuity stream. 
Obviously, $\wt Y_n$ has the same \ds\ as
\beao
Y_n&=&X_{11}+X_{21}X_{22}+X_{31}X_{32}X_{33}+\cdots +X_{n1}\cdots X_{nn}\,,\quad n\ge 1\,,
\eeao
and, under mild conditions,  the \seq\ $(Y_n)$ has the a.s. limit
\beam\label{eq:0}
Y= \sum_{n=1}^\infty \Pi_n\;\mbox{where $\Pi_n= \prod_{i=1}^n X_{ni}$ for $n\ge 1$.}
\eeam
We assume that the infinite series in \eqref{eq:0} converges a.s.
Since $\Pi_n\stas 0$ is a necessary condition for this convergence to hold we need that
\beao
\log |\Pi_n|= \sum_{i=1}^n \log |X_{ni}|\stas -\infty\,,\qquad \nto\,.
\eeao
Hence the random walk $(\log |\Pi_n|)$ has a negative drift, i.e., $\E [\log |X|]<0$, possibly infinite.
\par
Throughout this paper we assume that there exists a positive number $\alpha$ \st 
\beam\label{eq:cramer}
h(\alpha)=\E [|X|^\alpha]=1\,.
\eeam 
Assume for the moment, that $X\ge 0$ a.s.
By convexity of the \fct\ $h(s)$, $h(\alpha+\vep)>1$ and $h(\alpha-\vep)<1$ for small $\vep\in (0,\alpha)$ where we assume that
 $h(s)$ is finite in some neighborhood of $\alpha$. Then for positive $\vep$,
\beao
\E[|Y|^{\alpha+\vep}] \ge \E [|\Pi_n|^{\alpha+\vep}]= (h(\alpha+\vep))^n\,.
\eeao
The \rhs\ diverges to infinity as $\nto$, hence $\E[|Y|^{\alpha+\vep}]=\infty$.
We also have for $\alpha\le 1$ and $\vep\in (0,\alpha)$,  
\beao
\E [|Y|^{\alpha-\vep}]\le \sum_{i=1}^n \E [|\Pi_{n}|^{\alpha-\vep}]=\sum_{i=1}^n (h(\alpha-\vep))^n<\infty\,.  
\eeao
For $\alpha>1$ a similar argument with the Minkowski inequality shows that $\E [|Y|^{\alpha-\vep}]<\infty$.
\par
These observations on the moments indicate that $|Y|$ has some heavy tail in the sense that
certain moments are infinite. 
In this paper we will investigate the precise \asy\ behavior of $\P(\pm Y>x)$ as $\xto$.
It will turn out that, under \eqref{eq:cramer} and some additional mild assumptions,  
\beam\label{eq:precise}
\P(Y>x)\sim  \begin{cases}\dfrac 2 {m(\alpha)}\,\dfrac {\log x} {x^\alpha}\,,&\mbox{if }X\ge0\mbox{ a.s.,}\\
\dfrac{ 1}{m( \alpha)} \dfrac {\log x} {x^\alpha}\,,&\mbox{if }\P(X<0)>0,
\end{cases}
\quad \xto\,,\nonumber\\
\eeam
where $m(\alpha)= \E[|X|^\alpha \log |X|]$ is a positive constant.
In the case $\P(X<0)>0$ we also have $\P(Y>x)\sim \P(Y<-x)$ as $\xto$.
\par

\par
An inspection of \eqref{eq:0} shows that the structure of $Y$ is in a sense close to
\beao
Y'=1+ \sum_{n=1}^\infty \Pi'_n\;\mbox{where $\Pi'_n= \prod_{i=1}^n X_{i}$ for $n\ge 1$.}
\eeao
This structure has attracted a lot of attention; see the recent monograph 
Buraczewski et al. \cite{buraczewski:damek:mikosch:2016} and the references therein. 
Indeed, assuming $X$ and $Y'$ independent, it is easy to see that
the following fixed point equation holds:
\beam\label{eq:fix}
Y'\eqd X\,Y'+1\,.
\eeam
If this equation has a solution $Y'$ for given $X$ it is not difficult to see that
the stationary solution $(Y_t')$  
to the stochastic recurrence equation
\beam\label{eq:sre}
Y'_t= X_t\,Y_{t-1}'+1\,,\qquad t\in\bbz\,,
\eeam
satisfies \eqref{eq:fix} for $Y'=Y'_t$, and if $Y'$ solves \eqref{eq:fix} it has the stationary \ds\ of the 
\MC\ decribed in \eqref{eq:sre}.
\par
One of the fascinating properties of \eqref{eq:fix} and \eqref{eq:sre} is that, under condition \eqref{eq:cramer}, 
these equations generate power-law tail behavior. Indeed, if $X\ge 0$ a.s.
\beam\label{eq:14}
\P(Y'>x)\sim \dfrac{\E [(X\,Y'+1)^\alpha- (X\,Y')^\alpha]}{\alpha\,m(\alpha)}\,\dfrac 1 {x^\alpha},\qquad \xto\,,
\eeam
and if $\P(X<0)>0$ then
\beam\label{eq:14a}
\P(\pm Y'>x)\sim \dfrac{\E [|X\,Y'+1|^\alpha- |X\,Y'|^\alpha]}{2\alpha\,m(\alpha)}\,\dfrac 1 {x^\alpha},\qquad \xto\,,
\eeam
This follows from Kesten \cite{kesten:1973} who also proved \eqref{eq:14} and \eqref{eq:14a} for the linear combinations
of solutions to multivariate analogs of \eqref{eq:sre}. Goldie \cite{goldie:1991} gave an alternative proof 
of \eqref{eq:14} and \eqref{eq:14a}
and also derived the scaling constants for the tails.
\par
We will often make use of Kesten's  \cite{kesten:1973} Theorems A and B,   and Theorem 4.1 in Goldie \cite{goldie:1991}; cf. Theorem~2.4.4 and 2.4.7 in Buraczewski et al. \cite{buraczewski:damek:mikosch:2016}. For the reader's convenience,  we formulate these results here,
tailored for our particular setting. In the case $\P(X<0)>0$ we did not find a result of type \eqref{eq:3}
in the literature. Therefore we give an independent proof in Appendix~\ref{sec:appA}.
\bth\label{thm:kesten}
Assume the following conditions:
\begin{enumerate}
\item
The conditional law of $\log |X|$ given $\{X\ne 0\}$ is non-arithmetic.
\item
There exists $\alpha>0$ \st\ $\E|[X|^\alpha]=1$ and $\E|[X|^\alpha\log |X|]<\infty$.
\item
$\P(X\, x+1=x)<1$ for every $x\in\bbr$.
\end{enumerate}
If either $X\ge 0$ a.s. or $\P(X<0)>0$  hold then \eqref{eq:14} or \eqref{eq:14a} hold, respectively. In both cases,
there is a constant $c>0$ \st
\beam\label{eq:3}
\P\big(\max_{n\ge 1} \Pi_n' >x\big)\sim c\,x^{-\alpha}\,,\qquad \xto\,.
\eeam
\ethe
Here and in what follows, $c,c',\ldots$ stand for any positive constants whose values are not of interest.
\par
We have a corresponding result for the arithmetic case, i.e., when the law of 
$\log |X|$ conditioned on $\{X\ne 0\}$ is arithmetic.
This means that the support of $\log X$ (excluding zero if $\P(X=0)>0$) is a subset of $a\bbz$ for some non-zero $a$.
\bth\label{thm:grinc}
Assume conditions {\rm (2), (3)} of Theorem~\ref{thm:kesten} and 
\begin{enumerate}
\item[\rm (1')]
the law of $\log |X|$ conditioned on $\{X\neq 0\}$ is arithmetic.
\end{enumerate}
Then there exist constants $0<c<c'<\infty$ \st\ for large $x$,
\beam
x^\alpha \,\P\big(\max_{n\ge 1} \Pi_n' >x\big)&\in& [c,c']\,,\label{eq:3a}\\
x^\alpha\,\P\big(Y'' >x\big)&\in& [c,c']\,,\label{eq:3b}
\eeam
where $Y"= \sum_{n=1}^\infty|\Pi_n'|$
\ethe  
For $X\ge 0$, \eqref{eq:3a} is part of the folklore on ruin \pro y in the arithmetic case; see 
Asmussen \cite{asmussen:2003}, Remark~5.4, Section XIII. For the general case $\P(X<0)>0$ we refer
to the proof in Appendix~\ref{sec:appA}. 
Relation \eqref{eq:3b} can be found in 
Grincevi\v cius \cite{grincevicius:1975}, Theorem 2b.
\par
This paper has two main goals:
\begin{enumerate}
\item
We want to show that the \fct\ $\P(Y>x)$ is \regvary\ with index $-\alpha$ under the condition $\E [|X|^\alpha]=1$.
More precisely, we will show \eqref{eq:precise}.
\item
We want to show that 
\beam\label{eq:15}
\P(Y>x)\sim \sum_{n=1}^ \infty \P(\Pi_n>x)=:p(x)\,,\qquad x\to\infty\,.
\eeam
\end{enumerate}
Relation \eqref{eq:15} reminds one of similar results for sums of independent \regvary\ or subexponential
\rv s; see for example Chapter 2 in Embrechts et al. \cite{embrechts:klueppelberg:mikosch:1997}. The crucial
difference between \eqref{eq:15} and these results is that the summands $\Pi_n$ of $Y$ can be light-tailed
for every fixed $n$; the heavy tail of $Y$ builds up only for $\Pi_n$ with an index $n$ close to $\log x/m(\alpha)$.
\par
Positive solutions to these two problems are provided in Theorem~\ref{thm:1} and Corollary~\ref{cor:1}.
They also show that $\P(Y>x)/\P(Y'>x)\sim c \log x$ for some positive constant. The proof in Section~\ref{sec:proof} makes
use of Theorems~\ref{thm:kesten} and \ref{thm:grinc} as auxilary results.
We use classical exponential bounds for
sums of independent \rv s and change-of-\ms\ techniques; see Petrov's classic
\cite{petrov:1995} for an exposition of these results and techniques.
\par 
We also make an attempt to understand the tails of a vector-valued version of $Y$ when $\Pi_n=\bfX_{n1}\cdots\bfX_{nn}$ is the 
product of iid $d\times d$ matrices $(\bfX_{ni})$ with non-negative entries and a generic element $\bfX$ satisfies an analog of 
\eqref{eq:cramer} defining the value $\alpha>0$; see Section~\ref{subsec:kesten} for details. We define $\bfY=\bfY(\bfu)=\sum_{n=1}^\infty \Pi_n^\top\bfu$ for some unit vector $\bfu$ with 
non-negative components
and show that $\P(|\bfY|>x)$ is of the order $\log x/x^{\alpha}$. This approximation does not depend on the choice of $\bfu$ when $|\bfu|=1$. 
We prove this result
by showing the \asy\ equivalence between $\P(|\bfY|>x)$ and $p_\bfu(x)=\sum_{n=1}^\infty \P(|\Pi_n^\top\bfu|>x)$. Of course, the tail of $\bfY$ is
not characterized by the tail of the norm. Therefore we also consider linear combinations 
$\bfv^\top \bfY$ for any unit vector $\bfv$ with positve components and show that $\P(\bfv^\top \bfY>x)$ is also of the \asy\ order $\log x/x^\alpha$.
\par
This paper is structured as follows. In Section~\ref{sec:main} we present the main results  
 in the univariate case (Theorem~\ref{thm:1}
and Corollary~\ref{cor:1}) followed by a discussion of the results. 
Proofs are given in Section~\ref{sec:proof}. In Appendix~\ref{sec:appA} we provide proofs of
relations \eqref{eq:3} and \eqref{eq:3a} in the case when $\P(X<0)>0$; we did not find a corresponding result in the literature.
The multivariate case is treated in Section~\ref{sec:4}; Theorem~\ref{thm:3} is a multivariate analog of Theorem~\ref{thm:1} and Corollary~\ref{cor:1}.

\section{Main results}\label{sec:main}\setcounter{equation}{0}
We formulate one of the main results of this paper.
\bth\label{thm:1} Assume the conditions of Theorems~\ref{thm:kesten} or \ref{thm:grinc}, in particular there
exists $\alpha>0$ \st\ $h(\alpha)=\E [|X|^\alpha]=1$. In addition, we assume that $\E [|X|^\alpha (\log |X|)^2 ]<\infty$, or $\E [|X|^\alpha (\log |X|)^2 ]=\infty$
and $\E[|X|^\alpha \1(\log |X|>x)]$ is \regvary\ with index $\kappa\in (1,2]$, 
\begin{enumerate}
\item If $X\ge 0$ a.s.
then
\beao
p(x)\sim  \dfrac {2}{m(\alpha)} \dfrac{\log x}{x^\alpha}\,,\qquad \xto\,.
\eeao
\item If $\P(X<0)>0$ then
\beao
p(x)\sim  \dfrac {1}{m(\alpha)} \dfrac{\log x}{x^\alpha}\,,\qquad \xto\,.
\eeao
\end{enumerate}
\ethe
\bre\label{rem:1a}
In the course of the proof of Theorem~\ref{thm:1} we show that for $X\ge 0$ a.s. 
\beam\label{eq:18}
p(x)&\sim& \sum_{n=1}^{[\log x/m(\alpha)]} \P(\Pi_n>x)\nonumber\\& \sim& 2 x^{-\alpha}
\sum_{n=1}^{[\log x/m(\alpha)]}\Phi\big((\log x - n\,m(\alpha))/\sqrt{\sigma^2(\alpha) n}\big)\,,\qquad \xto\,,
\eeam
where $\Phi$ is the standard normal \df\ and $\sigma^2(\alpha)= \E[X^\alpha (\log X)^2]-(m(\alpha))^2$ is assumed finite.
\ere
\par
The following result is an immediate con\seq\ of Theorem~\ref{thm:1} and Proposition~\ref{prop:1}.
\bco\label{cor:1}
Assume the conditions of Theorem~\ref{thm:1}. If $X\ge 0$ a.s. then
\beam\label{eq:13}
\P(Y>x)\sim \sum_{n=1}^\infty \P(\Pi_n>x)\sim  \dfrac {2}{m(\alpha)} \dfrac{\log x}{x^\alpha}\,,\qquad \xto\,. 
\eeam 
If $\P(X<0)>0$ then
\beam\label{eq:13a}
\P(\pm Y>x)\sim \sum_{n=1}^\infty \P(\Pi_n>x)\sim  \dfrac {1}{m(\alpha)} \dfrac{\log x}{x^\alpha}\,,\qquad \xto\,. 
\eeam 
\eco
In contrast to the distinct results for $\P(Y'>x)$ in 
Theorems~\ref{thm:kesten} and \ref{thm:grinc}
for the non-arithmetic and arithmetic cases, respectively, relations \eqref{eq:13} and \eqref{eq:13a} hold in both
cases. In particular, in contrast to Theorem~\ref{thm:grinc} for $\P(Y'>x)$, we get precise \asy s for $\P(Y>x)$  in the arithmetic case.
Corollary~\ref{cor:1} and Kesten's Theorem~\ref{thm:kesten} in the general and in the non-arithmetic cases, respectively, show
that $\P(Y>x)$ and $\P(Y'>x)$ are
\regvary\ \fct s with index $-\alpha$. However, we have $\P(Y>x)/\P(Y'>x) \to \infty$ as $\xto$, accounting for the additional independence
of $(\Pi_n)$
in the structure of $Y$. In the non-arithmetic case 
we can even compare the scaling constants in the tails. For example, for $X\ge 0$ a.s. we have (see \eqref{eq:14})
\beao
\dfrac{\P(Y>x)}{\P(Y'>x)}&\sim& 
\dfrac{2\alpha}{ \E[ (X\, Y'+1)^\alpha-(X\,Y')^ \alpha]}\,\log x\,.
\eeao
\par
We proved \eqref{eq:13} under conditions implying that $\E [X^\alpha (\log X)^{1+\delta}]<\infty$ for some $\delta>0$ 
which is slightly stronger than the condition $m(\alpha)<\infty$ in Kesten's theorem.
\par
We observe the similarity of the results in Theorem~\ref{thm:kesten} and Corollary~\ref{cor:1} as regards the \asy\
symmetry of the tails in the case when $\P(X<0)>0$. In both cases, we have $\P(Y'>x)\sim \P(Y'<-x)$ and $\P(Y>x)\sim \P(Y<-x)$ as $\xto$.
Moreover, in this case we also have 
\beao
\P(|Y|>x)\sim \P\Big(\sum_{n=1}^\infty |\Pi_n| >x\Big)\,,\qquad \xto\,.
\eeao
\subsection{Implications and discussion of the results}
The tail behavior of $\P(Y>x)$ described by Corollary~\ref{cor:1} immediately ensures limit theory for the extremes and 
partial sums of an iid \seq\ $(Y_i)$ with generic element $Y$. Assuming the conditions of Theorem~\ref{thm:1} and $X\ge 0$ a.s., 
choose $a_n= (2n\,\log n/(\alpha\,m(\alpha)))^{1/\alpha}$. Then we know from classical theory that
\beam
a_n^{-1} \max_{i=1,\ldots,n} Y_i&\std& \xi_\alpha\,,\label{17a}\\
a_n^{-1} \Big(\sum_{i=1}^n Y_i-c_n\Big)&\std& S_\alpha\,;\label{17b}
\eeam
see for example Chapters 2 and 3 in Embrechts et al. \cite{embrechts:klueppelberg:mikosch:1997}.
Relation \eqref{17a} holds for any $\alpha>0$ and the \ds\ of $\xi_\alpha$ is Fr\'echet with parameter $\alpha$.
Relation \eqref{17b} holds only for $\alpha\in (0,2)$ and the \ds\ of $S_\alpha$ is $\alpha$-stable. The centering constants
$c_n$ can be chosen as $n\, \E[Y]$ for $\alpha>1$, $n\,\E[|Y|\1(Y\le a_n)]$ for $\alpha=1$ and $c_n=0$ for $\alpha\in (0,1)$.
\par
One can introduce the stationary \ts\ 
\beao
Y_n= \sum_{i=-\infty}^n \prod_{j=n-i+1}^n X_{ij}\,,\qquad n\in\bbz\,.
\eeao
We observe that $Y_n\eqd Y$. Unfortunately, $Y_n$ cannot be derived via an affine \sre\ as in the Kesten case;
see \eqref{eq:sre}. Therefore its dependence structure is less straightforward. However, it is another example of
a \ts\ whose power-law tails do not result from heavy-tailed input variables $X_{ni}$.
\par
Now assume for the sake of argument that $X\ge 0$ and $\log X$ has a non-arithmetic \ds . 
Write $S_n'=\log \Pi'_n=\sum_{i=1}^n \log X_i$. As a byproduct from
Theorem~\ref{thm:1} and \eqref{eq:3} we conclude that 
\beao
\sum_{n=1}^\infty \P\big(S_n'>x\mid \max_{j\ge 1} S_j' >x\big)=
\dfrac{\sum_{n=1}^\infty \P(S_n'>x)}{\P\big(\max_{j\ge 1} S_j' >x\big)}\sim c\,x\,.
\eeao
From \eqref{eq:18} and the latter relation we also obtain
\beao
\lefteqn{\dfrac 1 {x}\sum_{n=1}^{[x/m(\alpha)]} \P\big(S_n'>x\mid \max_{j\ge 1} S_j' >x\big)}\\
&=& \dfrac 1 x\,\E\Big[\#\{n\le [x/m(\alpha)]: S_n'>x \}\mid \max_{j\ge 1} S_j' >x\Big]
\to c\,,\qquad \xto\,.
\eeao
The constant $c$ can be calculated explicitly. Indeed, it has a  
nice interpretation in terms of a so-called {\em extremal index;} see Section~8.1 in Embrechts et al.
\cite{embrechts:klueppelberg:mikosch:1997} and Leadbetter et al. \cite{leadbetter:lindgren:rootzen:1983}
for its definition and properties. 
\par
Notice that the maxima of $(S'_t)$ have the same distribution as those of the {\em Lindley process}  given by
\beam\label{eq:stat}
S_t^{+}=\max(S_{t-1}^+ + \log X_t,0),\qquad t\ge 1,\qquad S_0^{+}=0\,;
\eeam
see Asmussen \cite{asmussen:2003}, Section III.6.
As $\E[\log X_0]<0$ the existence of the stationary solution $\tilde S_0^{+}$ to \eqref{eq:stat} 
is ensured since  $\{0\}$ is an atom. 
The extremal behavior of the Lindley process is well studied: its extremal index $\theta$ exists, is positive and satisfies 
\beao
 \E\Big[\#\{n\le [c'\log(x)]: \tilde S_n^+>x \}\mid \max_{1\le j\le [c'\log x]} \tilde S_j^+ >x\Big]
\to \frac1\theta\,,\qquad \xto\,,
\eeao
for some $c'>0$ depending on the exponential moments of the return time to the atom; see Rootz\'en \cite{rootzen:1988}. 
The extremal index can be expressed by using the Cram\'er constant for the associated ruin problem, i.e., 
the constant in \eqref{eq:3};  see Collamore and Vidyashankar \cite{collamore:2013}. From the previous discussion, we obtain
\beao
\dfrac 1 x\,\E\Big[\#\{n\le [x/m(\alpha)]: S_n^+>x \}\mid \max_{1\le j\le [x/m(\alpha)]} S_j^+ >x\Big]
\to \frac{2\alpha}{\theta}\,,\qquad \xto\,.
\eeao
\par
Surprisingly, under certain conditions the tail decay rate in \eqref{eq:13} is the same as for the solution to the  fixed point equation
\beao
\wt Y\eqd \sum_{i=1}^N X_i\,\wt Y_i\,,
\eeao
where $(\wt Y_i)$ are iid copies of $\wt Y$, $(X_i)$ is an iid positive \seq\ and $N$ is positive integer-valued. Moreover,
$(\wt Y_i)$, $(X_i)$ are mutually independent. In this case, the tail index $\alpha>0$ is given as the unique solution to the 
equation $\wt m(\alpha)=\E \big[\sum_{i=1}^N X_i^\alpha\big]=1$. The decay rate in \eqref{eq:13} is the same as for $\P(\wt Y>x)$ if $\alpha\in (0,1)$ and 
$\wt m'(\alpha)=0$. Results of this type appear in the context of smoothing transforms, branching and telecommunication models;
see Buraczewski et al. \cite{buraczewski:damek:mikosch:2016}, in particular Theorem~5.2.8(2), and the references therein.

\subsection{Examples}
\begin{figure}[htbp]
    \begin{center}
                \includegraphics[width=7cm]{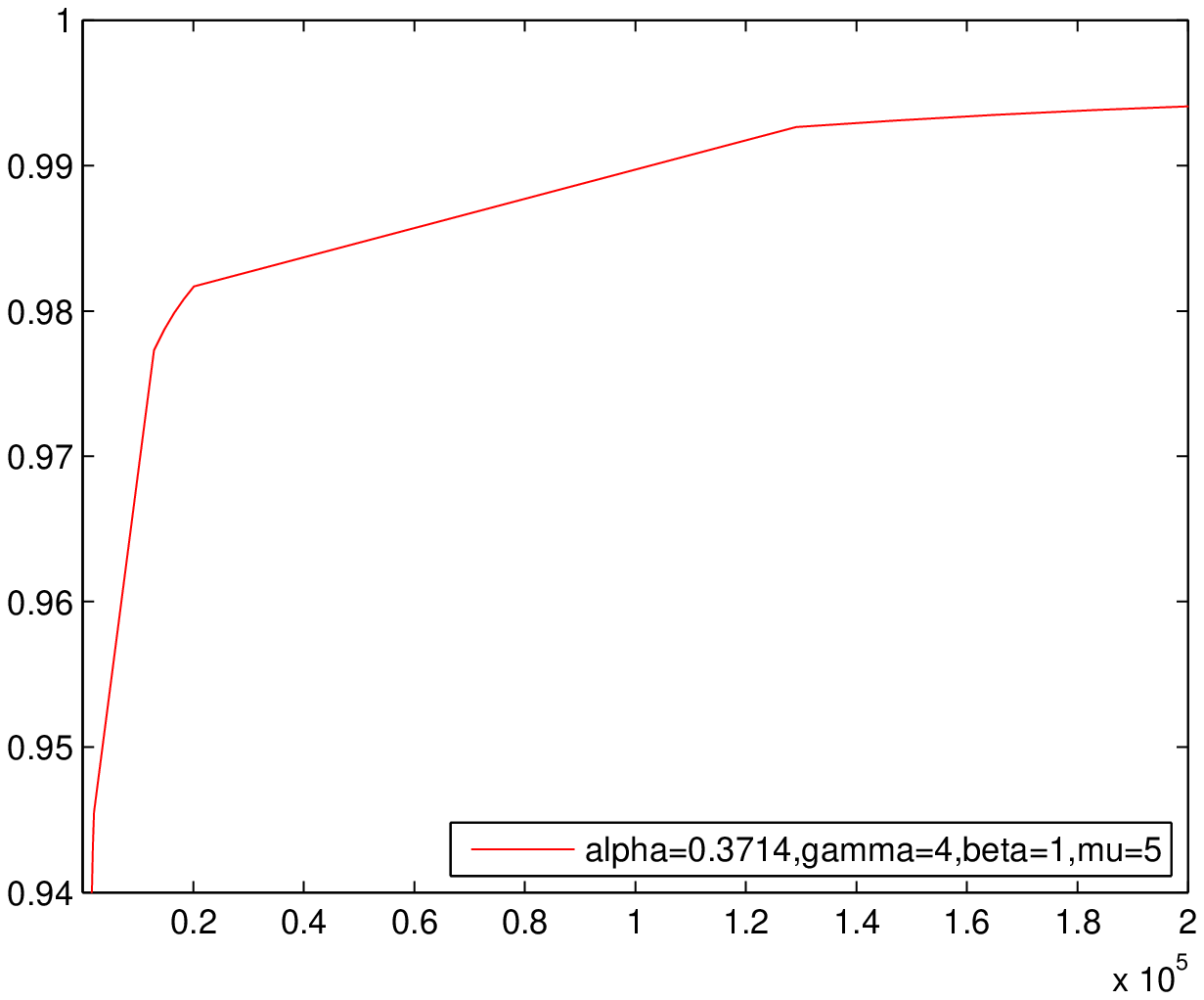}
        \includegraphics[width=7cm]{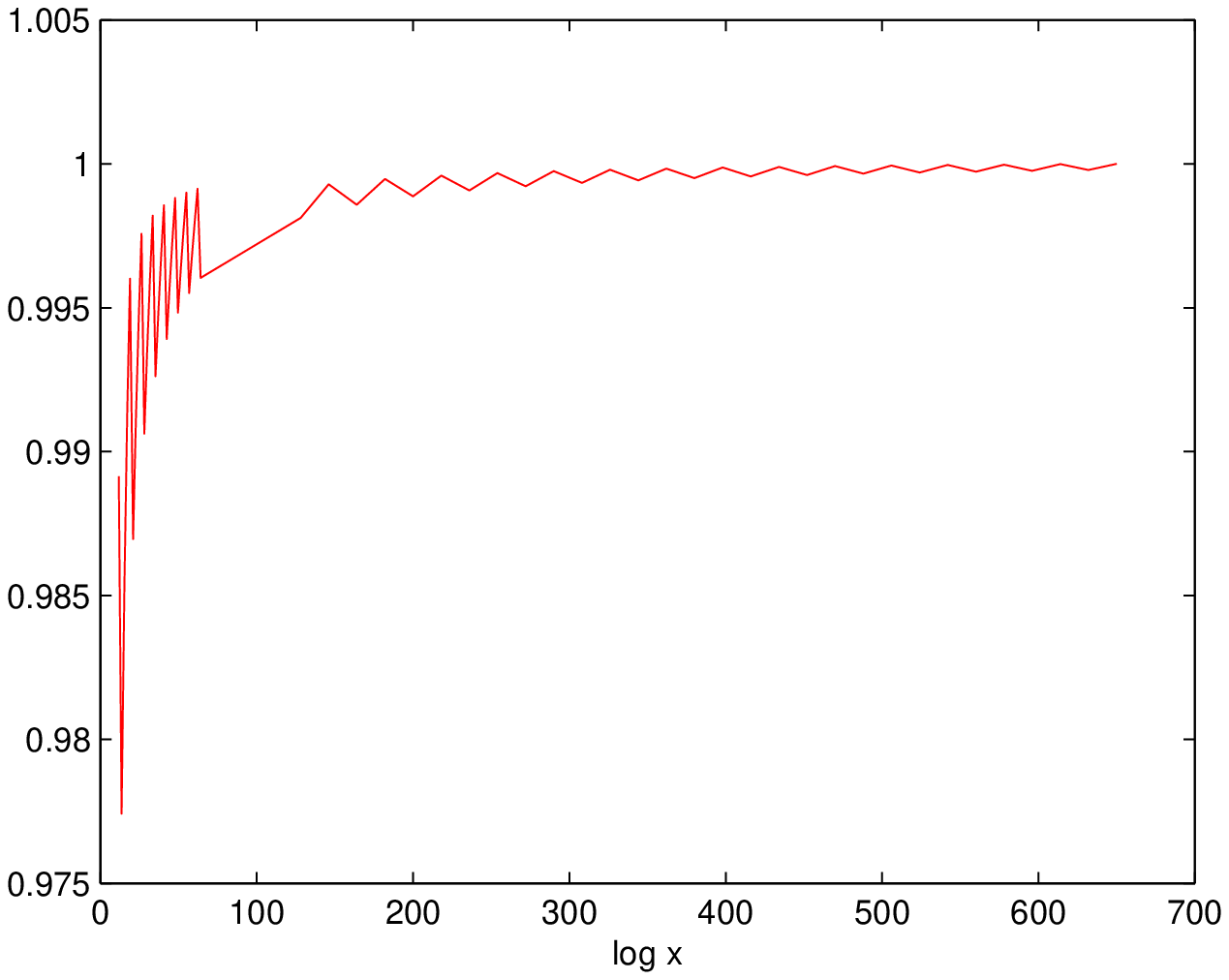}
    \end{center}
\caption{\small Two approximations of the ratio $x^\alpha p(x)/ (2\log x/m(\alpha))$ for loggamma-distributed $X$ with parameters
$\beta,\gamma,\mu$ and the resulting $\alpha$. The top figure shows the results of the normal approximation, the bottom figure 
a more precise approximation via change of \ms .}\label{Fig2a}
\end{figure}

In this section we illustrate our theory by considering various examples.  
\bexam\rm 
We assume that $(X_i)$ is an iid lognormal \seq , where $\log X$ has an $N(\mu,1)$ \ds\ with negative $\mu$.
Then for $s>0$,
\beao
\log (\E[X^s])=\mu\, s+s^2/2\,,\quad m(\alpha)=\alpha/2,,\quad \alpha=-2\mu,\quad \sigma^2(\alpha)=1\,.
\eeao
Notice that
$x^\alpha \,p(x)/2=\sum_{n=1}^{\infty} \Phi\big((\log x-n\,\alpha/2)/\sqrt{n}\big)$. 
\eexam
\bexam\rm
Assume that $X$  has  a $\Gamma(\gamma,\beta)$-density given by 
\beam\label{eq:gammadensity}
f_X(x)=\frac{\beta^\gamma x^{\gamma-1}\,\ex^{-x\beta}}{\Gamma(\gamma)}\,,\quad \beta,\gamma,x>0\,.
\eeam
Since $X$ has unbounded support the equation $\E[X^\alpha]=1$ always has a unique positive solution. For given values 
$\alpha$ and $\gamma$ we can determine suitable values $\beta$ \st 
\beao
\E[X^\alpha]=\frac{\Gamma(\gamma+\alpha)}{\Gamma(\alpha)\beta^\alpha}=1\,.
\eeao
We also have 
\beao
m(\alpha)&=&\frac{\Gamma(\gamma+\alpha)}{\beta^\alpha\Gamma(\gamma)}\Big(\dfrac{\Gamma'(\gamma+\alpha)}{\Gamma(\gamma+\alpha)}
-\log\beta\Big)\,,\\
\E[X^\alpha\,(\log X)^2]&=&\frac{\Gamma(\gamma+\alpha)}{\beta^{\alpha}\Gamma(\gamma)}
\Big( (\log\beta)^2-2\log\beta\dfrac{\Gamma'(\gamma+\alpha)}{\Gamma(\gamma+\alpha)}+
\dfrac{\Gamma''(\gamma+\alpha)}{\Gamma(\gamma+\alpha)}\Big)\,.
\eeao 
\eexam
\bexam\rm
Assume that $X=\ex^{Z-\mu}$ for some positive $\mu >0$ and a $\Gamma(\gamma,\beta)$-distributed \rv\ $Z$, i.e.,
$X$ has a loggamma \ds . For $\alpha<\beta$ we can calculate
\beao
\E[X^\alpha]=\ex^{-\alpha\mu}\,\big(1-\dfrac{\alpha}{\beta}\big)^{-\gamma}\,.
\eeao
The equation $\E[X^\alpha]=1$ has a positive solution \fif\ $\beta\mu>\gamma$.
Under this assumption,
\beao
m(\alpha)
&=&\ex^{-\alpha \mu}\big(1-\frac{\alpha}{\beta}\big)^{-\gamma-1}\dfrac 1\beta \big[\gamma+\mu(\alpha-\beta)\big]\,,\\
\sigma^2(\alpha)&=&\ex^{-\alpha \mu}\,\big(1-\frac{\alpha}{\beta}\big)^{-\gamma-2}\dfrac{1}{\beta^2}
\big[\gamma+(\gamma+\mu(\alpha-\beta))^2 \big]\,.
\eeao
Consider iid copies  $(Z_i)$ of $Z$. Then
\beao
p(x)=\sum_{n=1}^\infty \P\big(\sum_{j=1}^nZ_{j}-n\mu >\log x\big)\,,
\eeao
where $\sum_{j=1}^nZ_{j}$ is $\Gamma(n\gamma,\beta)$-distributed. In principle, this formula could be evaluated exactly by using the
gamma \df s. However, the events  $\{\sum_{j=1}^nZ_{j}>\log x\}$ are very rare for large $x$. Therefore one needs change-of-\ms\ techniques
to evaluate $p(x)$ or suitable approximation techniques. 
In the top Figure~\ref{Fig2a} we plot  the ratio of the normal approximation of $x^\alpha p(x)$ given in \eqref{eq:18} and $2\log x/m(\alpha)$
for $\mu=5$, $\gamma=4$ and $\beta=1$. For the same parameter set,
in the bottom figure we plot the ratio of 
\beao
x^{\alpha}p(x)=\sum_{n=1}^\infty \P^\alpha\Big(\sum_{j=1}^nZ_{j}>\log x+n\mu\Big)=\sum_{n=1}^\infty  F_{n}(\log x+n\mu)\,,
\eeao
where  we changed the \ms\ from $\P$ to $\P^\alpha(Z\in dx)= \ex^{\alpha x} \P(Z \in dx)$ resulting in the $\Gamma(n(\gamma+\alpha),\beta)$-\ds\ $F_n$.
The rationale for this change of \ms\ is explained in the proof of Theorem~\ref{thm:1}. A comparison 
of the two graphs shows the (not unexptected) result 
that the precise approximation of $x^\alpha p(x)$ via change of \ms\ is better than the approximation via the normal law.
\eexam

\section{Proof of Theorem~\ref{thm:1}}\label{sec:proof}\setcounter{equation}{0}
\subsection{First approximations}\setcounter{equation}{0}
Recall the definition of $p(x)$ from \eqref{eq:15}.
\bpr\label{prop:1} Assume the conditions of  Theorem~\ref{thm:1} and that $p(x)\sim c\,\log x /x^{\alpha}$. 
Then
\beam\label{eq:1a}
\P (Y>x)\sim p(x)\,,\qquad \xto\,. 
\eeam
\epr
A proof of the fact that $p(x)\sim c\,\log x /x^{\alpha}$ will be given in Section \ref{sec:precise}. 
\begin{proof}
Since $\sup_{i\ge 1}\P(\Pi_i>x)\le p(x)\to 0$ we have  
\beao
\P (Y>x)&\ge & \sum_{n=1}^\infty \P( \Pi_n>x)\,\P(\Pi_i\le x\,,i\ne n)\\
&=& \sum_{n=1}^\infty \P( \Pi_n>x) \prod_{i\ne n} \P(\Pi_i\le x)\\
&=&\sum_{n=1}^\infty \P( \Pi_n>x) \exp\Big(\sum_{i\ne n}^\infty \log (1- \P(\Pi_i> x) \Big)\\
&=&\sum_{n=1}^\infty \P( \Pi_n>x) \exp\Big(-(1+o(1))\sum_{i\ne n}^\infty\P(\Pi_i> x) \Big)\,.
\eeao
Hence
\beao
\P (Y>x)&\ge &p(x)\,\ex^{-(1+o(1)) p(x)}
=  (1+o(1))p(x)\,,
\eeao
and the liminf-part in \eqref{eq:1a} follows.
\par
Next we consider an upper bound for $\P(Y>x)$. We have for $\vep\in (0,1)$,
\beao
 \P(Y>x)&\le & p((1-\vep)x) + \P\big(C)\,.
\eeao
where $C=\{Y>x\,,\max_{n\ge 1} \Pi_n\le x(1-\vep)\}$.
We write for small $\delta>0$,
\beao
B_1&=& \bigcup_{1\le i<j} \{|\Pi_i|>\delta x,|\Pi_j|>\delta x\}\,,\\
B_2&=& \bigcup_{i=1}^\infty \{|\Pi_i|>\delta x\,,\max_{j\ne i} |\Pi_j|\le \delta x \}\,,\\
B_3&=& \{\max_{i\ge 1} |\Pi_i|\le \delta x \}\,.
\eeao
Observe that we have by Markov's inequality for $\gamma\in (0,\alpha)$,
\beam\label{eq:16new}
\wt p(x)=\sum_{n=1}^\infty \P(|\Pi_n|>x)\le x^{-\gamma} \sum_{n=1}^\infty (h(\gamma))^n\le c\,x^{-\gamma}\,.
\eeam
Hence for $\gamma\in (\alpha/2,\alpha)$ 
\beao
\dfrac{\P(C\cap B_1)}{p((1-\vep)x)}&\le& \dfrac{\big(\wt p(\delta x)\big)^2}{p((1-\vep)x)}
\le c\, x^{\alpha-2\gamma}\to 0\,,\quad \xto\,.
\eeao
Similarly, by independence of $Y-\Pi_n$ and $\Pi_n$ for any $n\ge 1$, and since \eqref{eq:16new} holds,
\beao
\dfrac{\P(C\cap B_2)}{p((1-\vep)x)}&\le& \dfrac{1}{p((1-\vep)x)}\sum_{n=1}^\infty \P(|Y-\Pi_n|>\vep x)\,\P(|\Pi_n|>\delta x)\\
&\le &\dfrac{(\wt p(\min(\vep,\delta) x))^2}{p((1-\vep)x)}
\le c\,x^{\alpha-2\gamma}\to 0\,,\qquad  \xto\,.
\eeao
We also have
\beao
\P(C\cap B_3)&\le& 
\P\Big( \sum_{n=1}^\infty |\Pi_n|>x\,,\max_{i} |\Pi_i|\le \delta x\big)\\&\le& \P\Big( \sum_{n=1}^\infty |\Pi_n|\,\1(|\Pi_n|\le \delta\,x)>x\Big)\,.
\eeao
Choose $\xi\in (0,1)$ and write $g(x)= [c_0\log x]$ for some positive constant $c_0$ to be chosen later. Then
\beao\P(C\cap B_3) &\le &
\P\Big(\sum_{n=1}^{g(x)} |\Pi_n|\1(|\Pi_n|\le \delta x)>x (1-\xi)\Big)\\
&&+\sum_{n=g(x)+1}^\infty \P\big(|\Pi_n|\1(|\Pi_n|\le \delta x)>x\, \xi^{n-g(x)}\,(1-\xi)\big)\\
&=&P_1(x)+P_2(x)\,.
\eeao
We have by Markov's inequality with $\gamma\in (\alpha/2,\alpha)$,
\beao
\dfrac{P_2(x)}{p((1-\vep)x)}&\le & c\,x^{\alpha-\gamma} 
\sum_{n=g(x)+1}^\infty \xi^{-(n-g(x))\gamma}\,\E [|\Pi_n|^{\gamma}]\\
&=& c\,x^{\alpha-\gamma} (h(\gamma))^{g(x)} \sum_{n=0}^\infty \big(\dfrac{h(\gamma)}{\xi^\gamma}\big)^n\\
&\le &c\,x^{\alpha-\gamma} (h(\gamma))^{g(x)} (1-\phi)^{-1}\,.
\eeao
Here we choose $\xi$ and $\gamma$ \st\ $\phi=h(\gamma)/\xi^\gamma<1$. Then the \rhs\ converges to zero if we 
choose $c_0>0$ sufficiently large.
\par
Next we find a bound for $P_1(x)$. We apply Prohorov's inequality; see Petrov \cite{petrov:1995}, p.~77. 
For this reason, we need bounds on the first and second moments of $S(x)=\sum_{n=1}^{g(x)} |\Pi_n|\1(|\Pi_n|\le \delta x)$.
\ble\label{lem:1} We have the following bounds
\beao
\E [S(x)]&\le& \left\{\barr{ll} h(1)\,(1-h(1))^{-1}<\infty\,,&\alpha>1\,,\\
 c\,(\log x)^2\,,&\alpha=1\,,\\
 c\,\log x\,x^{1-\alpha}\,,&\alpha\in (0,1)\,.
\earr\right.\\
\var(S(x))&\le& \left\{\barr{ll} h(2)(1-h(2))^{-1}<\infty\,,&\alpha>2\,,\\
 c\,(\log x)^2\,,&\alpha=2\,,\\c\,\log x\;x^{1-\alpha/2}\,,&\alpha\in (0,2)\,.
\earr\right.
\eeao
\ele
\begin{proof} We start with the bounds for $\E [S(x)]$.
(1) If $\alpha>1$, $\E [|\Pi_n|]= (h(1))^n<1$. Hence 
\beao
\E [S(x)]\le \sum_{n=1}^\infty  (h(1))^n= h(1) (1-h(1))^{-1}\,.
\eeao
(2) If $\alpha=1$ we use a domination argument.
Indeed, we have
\beao 
\E [|\Pi_n|\1(|\Pi_n|\le z)]|=\int_0^{z} \P(|\Pi_n|>y)\,dy\le \int_0^z\P(Y''>y)\,dy\,,\qquad z>0\,,
\eeao
where $Y"= \sum_{n=1}^\infty|\Pi_n'|$. By \eqref{eq:14}--\eqref{eq:14a} and \eqref{eq:3b}, respectively, we have $x\,\P(Y">x) \in [c,c']$ 
for constants $0<c<c'<\infty$ and large $x$. Hence
\beao
\E[S(x)]\le g(x) \E  [Y''\1(Y''\le \delta x)]\le c\,(\log x)^2\,.
\eeao
(3) A similar argument 
in the case $\alpha\in (0,1)$ shows that $x^\alpha \P(Y">x)\in [c,c']$ 
and $\E [Y"\1(Y"\le x)]\le c\,x^{1-\alpha}$ for large $x$. Hence $\E[S(x)]\le c\, x^{1-\alpha}\,\log x $.\\[2mm]
Our next goal is to find bounds for $\var(S(x))$. 
(1) If $\alpha>2$ then
\beao
\var(S(x))\le \sum_{n=1}^{g(x)} \E [\Pi_n^2]\le  h(2)(1-h(2))^{-1}<\infty\,.
\eeao
(2) Now assume $\alpha=2$. Then we have
\beao
\var(S(x))\le \sum_{n=1}^{g(x)} 
\E [\Pi_n^2\,\1(|\Pi_n|\le \delta x)]\,.
\eeao
The same domination argument as for $\E[S(x)]$ in the case $\alpha=1$ yields 
\beao
\E [\Pi_n^2\,\1(|\Pi_n|\le z)]&=& \int_0^z \P(\Pi_n^2>y) dy\\&\le& 
\int_0^z \P((Y'')^2>y)\,dy \\&=&\E [(Y'')^2\1((Y'')^2\le z)]\,,
\eeao
Hence 
\beao
\var(S(x))\le g(x)\,\E [(Y'')^2\1((Y'')^2\le \delta x)]\le c\,(\log x)^2\,.
\eeao
(3) Assume $\alpha\in (0,2)$. In this case, similar arguments as for $\alpha=2$ yield 
$\P((Y")^2>x)\le  c x^{-\alpha/2}$ and
 \beao
\var(S(x))\le c\,x^{1-\alpha/2}\log x\,.
\eeao
\end{proof}
\noindent
From Lemma~\ref{lem:1} we conclude that for large $x$,
\beao
P_1(x)&\le &\P\big(S(x)-\E [S(x)]>0.5\,x (1-\xi)\big)\,.
\eeao
\par
Now an application of Prohorov's inequality to the \rhs\ yields
\beao
P_1(x)&\le & \exp\Big( 
- \dfrac{0.5\,x (1-\xi)}{2(2\,\delta x)}{{\rm arsinh } 
\dfrac{(2\delta x)(0.5\,x (1-\xi))} {2 \var(S(x))}}\Big)\\
&=&\exp\Big( 
- \dfrac{1-\xi}{8\,\delta }{\rm arsinh } \big(0.5x^2 \,\delta\,(1-\xi)/\var(S(x))\big)\Big)\,,
\eeao
where $ \operatorname {arsinh} y=\log(y+{\sqrt {y^2+1}})\ge \log(2y)$ for positive $y$.

Now assume $\alpha>2$. Choose $\delta$ so small that $(1-\xi)/(8\,\delta)>\alpha$ and apply Lemma~\ref{lem:1} for large $x$,
\beao
{\rm arsinh } \big(0.5x^2 \,\delta/\var(S(x))\big)&\ge & 0.5 \log \big(x^2 \,\delta\,(1-\xi)(1-h(2))/h(2)\big)\\
&\sim & \log x\,,\qquad\xto\,,
\eeao
Hence we may conclude that 
\beam\label{eq:4}
\dfrac{P_1(x)}{p((1-\vep)x)}&\le & c\,x^\alpha\,P_1(x)\to 0\,,\qquad \xto\,.
\eeam
We proved the limsup-part in \eqref{eq:1a} for $\alpha>2$. For $\alpha=2$, using $\var(S(x))\le c\,(\log x)^2$, a slight modification of the 
Prohorov bound yields the same result. For $\alpha\in (0,2)$, we conclude from Lemma~\ref{lem:1} that for large $x$,
\beao
{\rm arsinh } \big(0.5x^2 \,\delta(1-\xi)/\var(S(x))\big)&\ge & 0.5 \log \big(x^{1-\alpha/2} \,\delta\,(1-\xi)\big)\\
&\sim & 0.5\,(1-\alpha/2)\,\log x\,,
\eeao
Now choose $\delta>0$ so small that  $(1-\alpha/2)(1-\xi)/(16\,\delta)>\alpha$ and then \eqref{eq:4}, hence \eqref{eq:1a} follows.
\end{proof}

\subsection{Preliminaries}\label{sec:neg}
Our next goal is to show that $p(x)\sim c \log x/x^\alpha$.
Since we will treat the cases $X\ge 0$ a.s. and $\P(X<0)>0$ in a similar way we will follow an idea of Goldie \cite{goldie:1991}. We define 
\beao
N_0=0\,,\quad N_i=\inf\{k>N_{i-1}: \Pi_k>0\}\,,\qquad i\ge 1\,\,,\quad I=\{N_i\,,i\ge 1\}\,,
\eeao
If $X\ge 0$ a.s. we have $N_i=i$ a.s. We introduce the non-negative variables $\tilde X_i=\prod_{j=N_{i-1}+1}^{N_{i}}|X_i|$, $i\ge 1$, and their products 
$\tilde \Pi_n=\prod_{i=1}^n \tilde X_i$  so that
\beao
p(x)&=&\sum_{n=1}^\infty \P(\Pi_n>x)= \E\Big[\sum_{n=1}^\infty \1(\Pi_n>x)\Big]\\
&=&
\E\Big[\sum_{n=1}^\infty \1(\tilde\Pi_n>x)\Big]=\sum_{n=1}^\infty \P(\tilde\Pi_n>x)\,.
\eeao
By independence of $(X_i)$, $(\tilde X_i)$ are iid as well. 
Under $\E[|X|^\alpha]=1$, the process $\big(\prod_{j= 1}^{t}|X_i|^\alpha\big)_{t\ge 1}$ is a martingale adapted to the 
filtration $\mathcal F_t=\sigma(X_i,i\le t)$. As $N_1$ is a stopping time \wrt\ this filtration we derive that $\E[\tilde X_1^\alpha]=1$ by an application of the stopping time theorem for martingales.\\
\par
Write $\tilde X$ for a generic element of $(\tilde X_i)$. 
We will use the following notation for $s>0$,  assuming these  moments are finite:
\beao
\tilde M(s)&=&\log \E[\tilde X^s]\,,\\
\tilde m(s)&=&\tilde  M'(s)=\dfrac{\E[|\tilde X|^s \,\log |\tilde X| ]}{\E[|\tilde X|^s]}\,,\\
\tilde \sigma^2(\alpha)&=&\tilde  m'(s)=\dfrac{ \E [|\tilde X|^s(\log |\tilde X|)^2]\,\E[|\tilde X|^s]- (\E [|\tilde X|^s\log |\tilde X|])^2}{(\E[|\tilde X|^s])^2}\,.
\eeao
The expression for the distribution of $\log\tilde X$ given $\{\tilde X>0\}$ can be derived by mimicing the arguments of Goldie 
\cite{goldie:1991}. Denote $p=\P(X>0)=1-q$ and
\beao
\gamma_\pm(dy)=\frac{\P (\pm X>0,\log |X|\in dy )}{\P (\pm X>0)}\,.
\eeao
Then we have 
\beam\label{eq:distr}
\P(\log \tilde X\in \cdot)=p\,\gamma_+(\cdot)+\sum_{n=2}^\infty q^2\, p^{n-2}\,\gamma_-^{(2)}(\cdot)\ast\gamma_+(\cdot)^{(n-2)}\,,
\eeam
where $\ast$ denotes the convolution operator.
\par
We introduce the tilted measure $\P^\alpha$:
\beao
d\P^{\alpha}(\log |X|\le y)=\ex^{\alpha y}d\P(\log |X|\le y)\,,
\eeao
and denote expectation and variance \wrt\ $\P^\alpha$ by $\E^\alpha$ and $\var^\alpha$, respectively.
Under the assumption $\E[|X|^\alpha (\log |X|)^2]<\infty$ the following moments are finite and positive
\beao
\E^\alpha[\log |X|]&=& \E[|X|^\alpha\log |X|]=m(\alpha)\,,\\
\var^\alpha(\log |X|)&=& \E^\alpha[(\log |X|)^2]-(\E^\alpha[\log |X|])^2\,,\\&=&
\E[|X|^\alpha\,(\log |X|)^2]-(m(\alpha))^2=\sigma^2(\alpha)\,.  
\eeao
\par
Since the tilted measure $\P^\alpha$ preserves the sum structure one has the identity
\beao
d\P^{\alpha}(\log |\tilde X|\le y)=\ex^{\alpha y}d\P(\log |\tilde X|\le y)\,.
\eeao
Then one can also check the existence of the quantities
\beao
\E^\alpha[\log |\tilde X|]&=& \E[|\tilde X|^\alpha\log |\tilde X|]=\tilde m(\alpha)\,,\\
\var^\alpha(\log |\tilde X|)&=& \E^\alpha[(\log |\tilde X|)^2]-(\E^\alpha[\log |\tilde X|])^2\,,\\&=&
\E[|\tilde X|^\alpha\,(\log |\tilde X|)^2]-(\tilde m(\alpha))^2=\tilde \sigma^2(\alpha)\,.
\eeao
These quantities coincide with $m(\alpha)$ and $\sigma^2(\alpha)$ when $X\ge 0$ a.s. Otherwise,  calculation yields
\beam\label{eq:tsig}
\nonumber
\left\{\barr{ll}\tilde m(\alpha)=&2m(\alpha)\,,\\ 
  \tilde \sigma^2(\alpha)=&2\sigma^2(\alpha)+\dfrac pq(2m(\alpha))^2\,.\earr\right.
\eeam
Finally, notice that $\P(N_1=1)=p$ and $\P(N_1=n)=q^2\,p^{n-2}$ so that $N_1$ admits finite moments of any order. Moreover, 
$\E[N_1]=1$ if $p=\P(X>0)=1$ and 2 else. If $\P^\alpha(\log |X|>x)$ is \regvary\ with index $\kappa\in (1,2]$, 
one can apply Case (b3) on p.~130 of Resnick \cite{resnick:1986} to the stopped random walk $\tilde X$ to obtain the equivalence
\beam\label{eq:forgot}
\P^\alpha(\log |\tilde X|>x)\sim \E[N_1]\,\P^\alpha(\log |  X|>x),\qquad \xto.
\eeam
Hence $\P^\alpha(\log |\tilde X|>x)$ is also \regvary\ with index $\kappa\in (1,2]$.

\subsection{More precise \asy s}\label{sec:precise}

On the set $\{\tilde \Pi_n>0\}$ we may write
\beao
\tilde S_n=\log \tilde \Pi_n= \sum_{i=1}^n \log \tilde X_i\,,\qquad n\ge 1\,.
\eeao
We will show that
\beao
p(x)=\sum_{n=1}^\infty \P(\tilde S_n-n\,\E [\log \tilde X]>\log x-n\,\E[\log \tilde X]\,,\tilde \Pi_n>0) \sim x^{-\alpha}\,L(x)  
\eeao
for a suitable \slvary\ \fct\ $L$.
Since $\E[\log \tilde X]=c\,\E[\log |X|]<0$, where $c=1$ for $X\ge 0$ and $c=2$ otherwise,
the random walk $(\tilde S_n)$ has negative drift. We will exploit the fact that, after the  change of \ms\
via $\P^\alpha$, the random walk $(\tilde S_n)$ has a positive drift $(n\,\tilde m(\alpha))$.
\par
In what follows, we will get bounds for sums of $\P(\tilde \Pi_n>x)$ over different $n$-regions. 
It will be convenient to use the following notation 
\beao
g_\xi(x)=[(1+\xi) \log x/\tilde m(\alpha)]\quad \mbox{for real $\xi$.}
\eeao
\ble\label{lem:gx} 
For any small $\vep>0$ we can find $\delta>0$ \st\ for sufficiently large  $x$,
\beam\label{eq:5a}
\sum_{n=g_\vep (x)}^\infty  \P(\tilde \Pi_n>x)\le c\,x^{-(\alpha+\delta)}\,.
\eeam
\ele
\begin{proof}
Denote $\tilde h(s)= \E[\tilde X^s]$ for $s\le \alpha$ and notice that $\tilde M(s)=\log\tilde h(s)$. By Markov's inequality for small $\vep\in (0,\alpha)$,
\beao\lefteqn{
x^\alpha\sum_{n=g_\vep (x)}^\infty \P(\tilde \Pi_n>x)}\\
&\le& x^{\vep}\,\sum_{n=g_\vep (x)}^\infty (\tilde h(\alpha-\vep))^n=x^{\vep}\,
(\tilde h(\alpha-\vep))^{g_\vep(x)} (1-\tilde h(\alpha-\vep))^{-1}\\
&=& \exp\Big(\log x\, \big(\vep  + \dfrac {[(1+\vep) \log x/\tilde m(\alpha)]}{\log x} \tilde M(\alpha-\vep)\big) \Big)(1-h(\alpha-\vep))^{-1}\,.
\eeao
By a  Taylor expansion, $\tilde M(\alpha-\vep)=\tilde M(\alpha-\vep)-\tilde M(\alpha)\sim -\tilde m(\alpha)\,\vep$ as $\vep\downarrow 0$. This proves  \eqref{eq:5a} for small $\vep$.
\end{proof} 
We apply \eqref{eq:5a} and the fact that $\sum_{n=g_0(x)+1}^{g_\vep(x)}\P(\tilde \Pi_n>x)\le c\,\vep \log x$ to show that
\beao
p(x)=\sum_{n=1}^\infty\P(\tilde \Pi_n>x)=\sum_{n=1}^{g_0(x)}\P(\tilde \Pi_n>x)
+ o(\log x)\,,\qquad \xto \,.
\eeao
Next define  
\beao
\nu(x)=\ex ^{\alpha \log x }\sum_{i=1}^{g_0(x)} \P(\tilde \Pi_n>x)=\ex ^{\alpha \log x }
\sum_{n=1}^{g_0(x)}\P(\tilde S_n >\log x).
\eeao
Then we have
\beam\label{eq:17}
\nu(x)&=&\ex ^{\alpha \log x}  \sum_{n=1}^{g_0(x)} \int_{\log x}^\infty d\P(\tilde S_n\le t)\nonumber\\
&=&   \int_{\log x}^\infty\ex^{-\alpha (t-\log x)}d\Big(\sum_{n=1}^{g_0(x)}  \P^{\alpha}(\tilde S_n\le t)\Big)\nonumber\\
&=& \int_0^ \infty \ex^{-\alpha  s } d\nu_{\alpha}(s+\log x)\,,
\eeam
where 
$
\nu_{\alpha}(y)= \sum_{n=1}^{g_0(x)}
\P^{\alpha}(\tilde S_n\le  y)\nonumber
$.
We want to show that
$\ell(x)=\nu_\alpha (\log x)$
is a \slvary\ \fct . More precisely, we want to show that $\ell (x)\sim  g_0(x)$. 
We have $ \ell(x)\le g_0(x)$. 
\par 
First assume  $\E^\alpha[(\log |X|)^2]=\E [(\log |X|)^{2} |X|^\alpha]<\infty$, then $\E^\alpha[(\log \tilde X)^2]<\infty$; see the expression \eqref{eq:forgot}. 
By the \clt\
under the \ms\ $\P^\alpha$
we have
\beao
 \ell(x)&=& \sum_{n=1}^{g_0(x)}  \P^{\alpha}\Big(\dfrac{\tilde S_n-n\,\tilde m(\alpha)}{\sqrt{n}\tilde \sigma(\alpha)}\le \dfrac{\log x-n\,\tilde m(\alpha)}
{\sqrt{n}\tilde \sigma(\alpha)} \Big)\nonumber\\
&=&g_0(x)-\sum_{n=1}^{g_0(x)} \ov\Phi\Big(\dfrac{\log x-n\,\tilde m(\alpha)}{\sqrt{n}\tilde \sigma(\alpha)} \Big)+o(\log x)\\&=:& g_0(x)-T(x)+o(\log x)\,,
\eeao
where $\ov \Phi=1-\Phi$ denotes the right tail of the standard normal \ds .  We have
\beao
T(x)
&=&O(K_0)+\sum_{n=1}^{g_0(x)-K_0} \ov \Phi\Big(\dfrac{\log x-n\,\tilde m(\alpha)}{\sqrt{n}\tilde \sigma(\alpha)} \Big)\\
&\le &O(K_0)+(g_0(x)-K_0)\ov \Phi (K)
\,,
\eeao
where for a given $K>0$ we choose an integer $K_0$ so large that  $(\log x-n\,\tilde m(\alpha))/(\sqrt{n}\tilde\sigma(\alpha))>K$.
Since we can choose $K$ as large as we wish we finally proved that $\ell(x)\sim g_0(x)$.
\par 
Now assume that $\P ^\alpha(\log |X|>x)$ is \regvary\ for some $\kappa\in (1,2]$ and if $\kappa=2$ also assume 
that $\E^\alpha [(\log |X|)^2]=\infty$. Then we also have that $\P ^\alpha(\log \tilde X>x)$ is \regvary\ for some $\kappa\in (1,2]$ and if $\kappa=2$  
that $\E^\alpha [(\log \tilde X)^2]=\infty$; see the discussion at the end of Section \ref{sec:neg}. Choose $(a_n)$ \st\ 
\beao
n\,\big[
\P^\alpha (\log \tilde X>a_n)+a_n^{-2}\E^\alpha [(\log \tilde X)^2\1(\log \tilde X\le a_n)]\big]=1\,.
\eeao 
Then under $\P^\alpha$, 
\beao
a_n^{-1} (\tilde S_n-n\,\tilde m(\alpha))\std S_\kappa\,,
\eeao
where $S_\kappa$ has a $\kappa$-stable \ds\ $\Phi_\kappa$. The same arguments as above yield
\beao
 \ell(x)&=& \sum_{n=1}^{g_0(x)}  \P^{\alpha}\Big(\dfrac{\tilde S_n-n\,\tilde m(\alpha)}{a_n}\le \dfrac{\log x-n\,\tilde m(\alpha)}
{a_n} \Big)\nonumber\\
&=&g_0(x)-\sum_{n=1}^{g_0(x)} \ov \Phi_\kappa\Big(\dfrac{\log x-n\,\tilde m(\alpha)}{a_n} \Big)+o(\log x)\\&=:& g_0(x)-T(x)+o(\log x)\,.
\eeao
and also $\ell(x)\sim g_0(x)$.
\par
Finally, integrating by parts, we obtain from \eqref{eq:17},
\beao
\nu(x)&=& \int_0^\infty \ex^{-\alpha  t } d\ell(x\ex^t)= \ell( x )+\alpha  \int_0^\infty \ex^{-\alpha  t }  \ell(x\ex^t)dt\\
&\sim & g_0( x )+\alpha  \int_0^\infty \ex^{-\alpha  t }  g_0(x\ex^t)\,dt\\&\sim&  \dfrac{2\log x} {\tilde m(\alpha)}
=\left\{\barr{ll}  \dfrac{2}{m(\alpha)}\,\log x &\mbox{if $X\ge 0$}\,,\\
\dfrac{1}{m(\alpha)}\,\log x &\mbox{otherwise}\,.\earr\right.
\eeao
This finishes the proof of the desired bounds for $p(x)$ and concludes the proof of Theorem~\ref{thm:1}.

\section{The multivariate case}\label{sec:4}\setcounter{equation}{0}
\subsection{Kesten's multivariate setting}\label{subsec:kesten}
In this section we will work under the conditions of the multivariate 
setting of Kesten \cite{kesten:1973}; see Section 4.4 in Buraczewski et al. \cite{buraczewski:damek:mikosch:2016}.
\par
We consider iid $d\times d$ matrices $(\bfX_t)$ with a generic element $\bfX$ \st\ $\P(\bfX\ge {\bf0})=1$ and $\bfX$ does not have a zero
row with \pro y~1. Here and in what follows, vector inequalities like $\bfx\ge \bfy$, $\bfx>\bfy,\ldots$, in $\bbr^d$ are understood 
componentwise. 
We write  $\bbs^{d-1}=\{\bfx\in\bbr^d: |\bfx|=1\}$ for the unit sphere in $\bbr^d$
and $\bbs_+^{d-1}=\{\bfx\in\bbr^d: |\bfx|=1\,,\bfx\ge \bf0\}$. We always use the Euclidean norm $|\cdot|$ and write 
$\|\cdot\|$ for the corresponding operator norm. 
Here and in what follows, all vectors are column vectors.
We write $\Pi_n'=\bfX_1\cdots \bfX_n$, $n\ge 1$.
\par
Assume the following conditions.
\begin{enumerate}
\item The top Lyapunov exponent $\gamma$ is negative:
\beao
\gamma= \inf_{n\ge 1} n^{-1} \E [\log\|\Pi_n'\|]<0\,.
\eeao
\item
Consider 
\beao
h(s)= \inf_{n\ge 1} \big(\E [\|\Pi_n'\|^s]\big)^{1/n}=\lim _{\nto} \big(\E [\|\Pi_n'\|^s]\big)^{1/n}
\eeao
and assume that there is $\alpha>0$ \st\ $h(\alpha)=1$. 
\item
$\E\big[\|\bfX\|^\alpha \log^+ \|\bfX\|\big]<\infty$.
\item
The additive subgroup of $\bbr$ generated by the numbers $\log \la_\bfs$ is dense in $\bbr$, where
$\la_\bfs$ is the dominant eigenvalue of $\bfs=\bfx_1\cdots \bfx_n$ for $\bfx_i$, $i=1,\ldots,n$, for some $n\ge 1$,  in the support of $\bfX$ and \st\ $\bfs>\bf0$.
\end{enumerate}
Let the $\bbr^d$-dimensional column vector $\bfY'$ be independent of $\bfX$.
Under the conditions above, the fixed point equation $\bfY'\eqd \bfX \bfY'+\bfu$ has
a unique solution $\bfY'$, where $\bfu\in \bbs^{d-1}$ is a deterministic vector. Then one has the \rep
\beao
\bfY'^\top\eqd\bfu^\top\sum_{n=1}^\infty \Pi_{n-1}'=\bfu^\top+\bfu^\top\sum_{n=2}^\infty \Pi_{n-1}'\,.
\eeao
\par
The next result follows from Theorems 4 and A in Kesten \cite{kesten:1973}; cf. Theorem~4.4.5 in Buraczewski et al.
\cite{buraczewski:damek:mikosch:2016}, where the conditions below are also compared with those in the original paper.
\bth\label{thm:kesten2}
Under the conditions above, there exists a finite \fct\ $e_\alpha$ on $\bbs^{d-1}$ \st
\beam\label{eq:kestena}
\lim_{\xto}x^\alpha \,\P(\bfv^\top\bfY' >x) =e_\alpha(\bfv)\,,\qquad \bfv\in \bbs^{d-1}\,,
\eeam
and $e_\alpha(\bfv)>0$ for $\bfv\in \bbs_+^{d-1}$.
Moreover, the limits
\beam\label{eq:kest3}
\lim_{\xto}x^\alpha\,\P\big(\max_{n\ge 1} |{\Pi_n'}^\top \bfu|>x\big)=\wt e_\alpha(\bfu)\,,\qquad \bfu\in \bbs^{d-1}\,,
\eeam
exist, are finite and positive for $\bfu\in \bbs_+^{d-1}$.
\ethe
We notice that this result is analogous to Theorem~\ref{thm:kesten} in the non-arithmetic case.
Indeed, it is a special case for $d=1$.
\bre\label{rem:kest} Under the condition of non-negativity on $\bfX$ and $\bfu$, \eqref{eq:kestena} implies \regvar\ of $\bfY'$
(see Buraczewski et al. \cite{buraczewski:damek:mikosch:2016}, Theorem~C.2.1),
i.e., there exists a Radon \ms\ $\nu$ on $\ov \bbr^d_0= (\bbr\cup \{\pm  \infty\})^d\backslash \{\bf0\}$ \st
\beao
x^\alpha\,\P(x^{-1} \bfY'\in \cdot)\stv \nu(\cdot)\,.
\eeao
Here $\stv$ denotes vague \con\ in $\ov \bbr^d_0$
and $\nu$ has the property $\nu(t\cdot)=t^{-\alpha}\nu(\cdot)$, $t>0$. In particular, for fixed $\bfu$ and any $\bfv\in\bbs^{d-1}_+$ there exist positive constants 
$c_\bfu$ and $c_{\bfu,\bfv}$ \st\ as $\xto$,
\beao
x^\alpha\,\P(|\bfY'|>x) \to c_\bfu\quad\mbox{and}\quad x^\alpha\,\P(\bfv^\top \bfY'>x) \to c_{\bfu,\bfv}\,.
\eeao
Under non-negativity of $\bfX$, $\bfu$ and $\bfv$, and if $\blue e_\alpha(\bfv)\ne 0$ for some $\bfv$,
\eqref{eq:kestena} still implies \regvar\ of $\bfY'$ for non-integer-valued $\alpha$; see 
 the comments after Theorem~C.2.1 in \cite{buraczewski:damek:mikosch:2016}.
\ere
\subsection{Main results}
In what follows, we provide an analog of the univariate theory built in the previous sections. For this reason,
consider an iid array $(\bfX_{ni})_{n,i=1,2,\ldots}$ with generic element $\bfX$. Assume the conditions on $\bfX$ and $\bfu$ from
Kesten's Theorem~\ref{thm:kesten2} and define
\beao
\bfY^\top= \bfY^\top(\bfu)=\bfu^\top\sum_{n=1}^\infty \Pi_{n}\,,\quad\mbox{where $\Pi_n= \prod_{j=1}^n \bfX_{nj}$.}
\eeao
For any unit vectors $\bfu,\bfv\in \bbs_+^{d-1}$, we define
\beao
p_\bfu(x)= \sum_{n=1}^\infty \P(|\Pi_n^\top \bfu|>x)\quad\mbox{and}\quad
p_{\bfu,\bfv}(x)= \sum_{n=1}^\infty \P(\bfv^\top \Pi_n^\top \bfu>x)\,.
\eeao
The following result is an analog of Theorem~\ref{thm:1}.  
\bth\label{thm:3} Assume the Kesten conditions of Section~\ref{subsec:kesten}, in particular there
exists $\alpha>0$ \st\ $h(\alpha)=1$. In addition, we assume that $\E [\|\bfX\|^\alpha (\log \|\bfX\|)^2 ]<\infty$.
Then we have for $\bfY=\bfY(\bfu)$ and any $\bfu,\bfv\in \bbs_+^{d-1}$ 
\beam\label{eq:wwa}
\P(|\bfY|>x)&\sim &p_\bfu(x)\sim  \dfrac {2}{m(\alpha)} \dfrac{\log x}{x^\alpha}\,,\qquad \xto\,,\\
\label{eq:ww}
\P( \bfv^\top \bfY>x) &\sim&  p_{\bfu,\bfv}(x) \sim  \dfrac {2(\bfv^\top\bfu)^\alpha}{m(\alpha)} \dfrac{\log x}{x^\alpha}\,,\qquad \mbox{if also}\; \bfu>\bf0,
\eeam
where $m(\alpha)=h'(\alpha)$ is independent of $\bfu,\bfv$.
\ethe
The proof of this result is given in Section~\ref{subsec:proofthm3}.
\bre
As in the univariate case it is possible to relax the condition $\E [\|\bfX\|^\alpha (\log \|\bfX\|)^2 ]<\infty$
by a \regvar\ condition of order $\kappa\in (1,2]$, assuming for $\kappa=2$ that  $\E [\|\bfX\|^\alpha (\log \|\bfX\|)^2 ]=\infty$.
This \regvar\ condition has to be required under the \pro y \ms\ $\P^\alpha$ which will be explained in the course of the
proof of the theorem. Write $Z=\sum_{ij} X_{ij}$ and $V=\min_{i=1,\ldots,d} \sum_{j=1} ^dX_{ij}$, where $X_{ij}$ are the entries of $\bfX$.
Then one needs to assume that as $\xto$,
\beam\label{eq:regvar}
\P^\alpha(\pm \log Z>x)&\sim& c_\pm\,\dfrac{L(x)}{x^\kappa}\;\; \mbox{and}\;\; \P^\alpha( \log V \le -x)=O(\P^\alpha(|\log Z|>x))\,,\nonumber\\~
\eeam
where $c_\pm$ are non-negative constants \st\ $c_++c_-=1$ and $L$ is a \slvary\ \fct . In the case when  $\E [\|\bfX\|^\alpha (\log \|\bfX\|)^2 ]<\infty$
we use a \clt\ with Gaussian limit of Hennion \cite{hennion:1997}. Under \eqref{eq:regvar} and $\P^\alpha$, 
one can instead apply a corresponding result with a $\kappa$-stable
limit. 
The corresponding results can be found in Hennion and Herv\'e \cite{hennion:herve:2008}; see their Theorem 1.1 (replacing 
Theorem~3 in  \cite{hennion:1997}) and Lemma~2.1 (replacing Lemma~5.1 in \cite{hennion:1997}). Thus, as in the univariate
case, the moment condition $\E [\|\bfX\|^\alpha (\log \|\bfX\|)^2 ]<\infty$ can be slightly relaxed. 
\par
In \cite{hennion:1997} and \cite{hennion:herve:2008} the condition
\beam\label{eq:hennion}
\P\big(\Pi_n>{\bf0}\;\mbox{for some $n\ge 1$}\big)>0\,.
\eeam
is assumed. This condition follows under the conditions of Section~\ref{subsec:kesten}; see p.~171 in \cite{buraczewski:damek:mikosch:2016}.
\ere
\bre We observe that 
\beao
\dfrac{p_{\bfu,\bfv}(x)}{p_{\bfu}(x)}\to (\bfv^\top \bfu)^\alpha\,,\qquad\xto\,.
\eeao
The \rhs\ is smaller than one unless $\bfu=\bfv$. In particular, for $\bfu=\bfv>\bf0$ we have
\beao
\P(|\bfY(\bfu)|>x) &\sim& 
\P(\bfu^\top \bfY(\bfu)>x)\,,\\
\P(|\bfY(\bfu)|>x, \bfu^\top \bfY(\bfu)\le x)&=& o(\log x/x^\alpha)\,.
\eeao
This means that $\bfY(\bfu)$ puts most tail mass in the direction of $\bfu$.
\ere
Following Remark~\ref{rem:kest}, we also have full \regvar\ of the vector $\bfY(\bfu)$ if $\bfu>0$
since  \eqref{eq:ww} holds for $\bfv\ge \bf0$.
\bco
Assume the conditions of Theorem~\ref{thm:3} and that $\alpha$ is not an integer. 
Then  $\bfY(\bfu)$ is \regvary\ with index $\alpha$. In particular, there is a 
Radon \ms\ $\nu$ on $\ov\bbr_0^d$ \st\
\beao
\dfrac{m(\alpha)}{2}\dfrac{x^\alpha}{\log x} \P(x^{-1}\bfY(\bfu)\in \cdot)\stv \nu(\cdot)\,,\quad \xto\,,
\eeao
and $\nu$ is uniquely determined by its values on the sets $A_\bfv=\{\bfy\in \bbr^d: \bfv^\top \bfy\}$, $\bfv\ge 0$, i.e., $\nu(A_\bfv)= (\bfv^\top \bfu)^\alpha$.
\eco

\subsection{Proof of Theorem~\ref{thm:3}}\label{subsec:proofthm3} The proofs of \eqref{eq:wwa} and \eqref{eq:ww} are
very much alike. We focus on \eqref{eq:wwa} and only indicate the differences with the proof of \eqref{eq:ww}. 
We will follow the lines of the 
proof of Theorem~\ref{thm:1} and Corollary~\ref{cor:1}.
\par
With start with an analog of Proposition~\ref{prop:1}.
\bpr\label{prop:2}
Assume the conditions of Section~\ref{subsec:kesten} and that $p_\bfu$, $p_{\bfu,\bfv}$ are \regvary .
Then for any $\bfu,\bfv\in \bbs^{d-1}_+$ and $\bfY=\bfY(\bfu)$,
\beam\label{eq:kest4}
\dfrac{\P(|\bfY|>x)}{p_\bfu(x)}\sim \dfrac{\P(\bfv^\top \bfY>x)}{p_{\bfu,\bfv}(x)}\sim 1\,,\qquad \xto\,.
\eeam
\epr
\begin{proof}
We follow the lines of the proof of Proposition~\ref{prop:1}.
Since $h(\alpha)=1$, we have by convexity of $h$ for $\gamma\in (0,\alpha)$, $h(\gamma)<1$, hence for sufficiently large $n$,
\beam\label{eq:expon}
\E [\|\Pi_n\|^\gamma]<c_0^n\,,
\eeam
for some  $c_0\in (0,1)$. By Markov's inequality, with $\wt p(x)= \sum_{n=1}^\infty \P(\|\Pi_n\|>x)$,
\beam\label{eq:kestenb}
x^\gamma p_{\bfv,\bfu}(x)\le x^\gamma p_\bfu(x)  &\le& x^\gamma\,\sum_{n=1}^\infty \P( \|\Pi_n\|>x)=:\wt p(x)\le  c\,.
\eeam
In particular, $p_\bfu(x)\to 0$. By \eqref{eq:kest3}
we have 
\beao
\sup _{i\ge 1}\P(|\Pi_i^\top\bfu|>x)\le \wt p(x)\to 0\,,
\eeao
Using this fact and the same arguments as in the univariate case, we derive
\beao
\P(|\bfY|>x)&\ge& \sum_{n=1}^\infty \P(|\Pi_n^\top \bfu|>x)\,\P(|\Pi_i^\top \bfu|\le x\,,i\ne n)\\
&=&p_\bfu(x) (1+o(1))\,,\\
\P(\bfv^\top \bfY>x)&\ge& p_{\bfu,\bfv}(x) (1+o(1))\,.
\eeao
This proves the liminf-part of \eqref{eq:kest4}.
\par
Next we prove the limsup-part.
We have for $\vep\in (0,1)$,
\beao
 \P(|\bfY|>x)&\le & p_\bfu((1-\vep)x) + \P\big(C)\,,
\eeao
where $C=\{|\bfY|>x\,,\max_{n\ge 1} |\Pi_n^\top \bfu|\le x(1-\vep)\}$.
We write for small $\delta>0$,
\beao
B_1&=& \bigcup_{1\le i<j} \{|\Pi_i^\top\bfu|>\delta x,|\Pi_j^\top\bfu|>\delta x\}\,,\\
B_2&=& \bigcup_{i=1}^\infty \{|\Pi_i^\top\bfu|>\delta x\,,\max_{j\ne i} |\Pi_j^\top\bfu|\le \delta x \}\,,\\
B_3&=& \{\max_{i} |\Pi_i^\top\bfu|\le \delta x \}\,.
\eeao
By \eqref{eq:kestenb} we have for $\gamma\in (0,\alpha)$ and large $x$, $ p_\bfu(x)\le x^{-\gamma}$.
Therefore we may proceed as in the univariate case and obtain for $\gamma\in (\alpha/2,\alpha)$,
\beao
\dfrac{\P(C\cap B_1)}{p_\bfu((1-\vep)x)} \le c\,x^{\alpha-2\gamma}\to 0\,,\quad \xto\,.
\eeao
Similarly,
\beao
\dfrac{\P(C\cap B_2)}{p_\bfu((1-\vep)x)}\le 
\sum_{n=1}^\infty \dfrac{\P(|\bfY-\Pi_n^\top\bfu|>\vep x)\,\P(|\Pi_n^\top \bfu|>\delta\,x)}{p_\bfu((1-\vep)x)}\le c x^{\alpha-2\gamma}\to 0\,,
\eeao
and for $\xi\in (0,1)$ and $g(x)=[c_0\log x]$, $c_0>0$,
\beao\P(C\cap B_3) &\le &
\P\Big(\sum_{n=1}^{g(x)} |\Pi_n^\top\bfu|\1(|\Pi_n^\top\bfu|\le \delta x)>x (1-\xi)\Big)\\
&&+\sum_{n=g(x)+1}^\infty \P\big(|\Pi_n^\top\bfu|\1(|\Pi_n^\top\bfu|\le \delta x)>x\, \xi^{n-g(x)}\,(1-\xi)\big)\\
&=&P_1(x)+P_2(x)\,.
\eeao
The proof of $P_2(x)/p_\bfu((1-\vep)x) \to 0$ is analogous to the univariate case. Write
\beao
S(x)= \sum_{n=1}^{g(x)} |\Pi_n^\top\bfu| \,\1(|\Pi_n^\top\bfu|\le \delta x).
\eeao
We have  similar bounds for $\E [S(x)]$ and $\var(S(x))$ as in the univariate case; see Lemma~\ref{lem:1}. The key to
this fact is domination via the inequalities
\beao
\E\big[ |\Pi_n^\top\bfu| \,\1(|\Pi_n^\top\bfu|\le z)\big]&=&\int_0^z\P(|\Pi_n^\top\bfu|>y)\,dy\\&\le&
\int_0^z\P\big(\max_{i\ge 1} |{\Pi'_i}^\top\bfu|>y\big)\,dy\,,\\
\E [|\Pi_n^\top\bfu|^2\1(|\Pi_n^\top\bfu|\le z)]&=& \int_0^z \P(|\Pi_n^\top\bfu|>\sqrt{y})\,dy\\
&\le& \int_0^z \P\big(\max_{i\ge 1} |{\Pi'_i}^\top\bfu|>\sqrt{y}\big)\,dy\,.
\eeao
Now exploit the result for the tails in  \eqref{eq:kest3} and the same domination argument as in the univariate case. 
Finally, Prohorov's inequality  applies to show that $P_1(x)/p_\bfu((1-\vep)x)\to 0$. 
\par
The proof of $\limsup_{\xto} \P(\bfv^\top \bfY>x)/ p_{\bfu,\bfv}((1-\vep)x)\le 1$ is analogous.
\end{proof}
\par
Our next goal is to show that $p_\bfu$ and $p_{\bfu,\bfv}$ are \regvary\ \fct s.
For $\vep >0$, we write $g_\vep(x)= [(1+\vep)\log x/m(\alpha)]$  where $m(\alpha)= h'(\alpha)$.
We continue with the following analog of Lemma~\ref{lem:gx}.
\ble\label{lem:4}
Assume the conditions of Section~\ref{subsec:kesten}.
Then
\beao
x^\alpha\sum_{n=g_\vep(x)}^\infty \P( |\Pi_n^\top\bfu|>x)= o(\log x)\,,\qquad \xto\,.
\eeao 
\ele
\begin{proof}
By Markov's inequality and a Taylor expansion of $\log h(\gamma)-\log h(\alpha)$ at $\alpha$ for $\gamma<\alpha$ close to $ \alpha$, since $\log h(\gamma)<0$,
\beao
x^\alpha\sum_{n=g_\vep(x)}^\infty \P( |\Pi_n^\top\bfu|>x)
&\le  & x^{\alpha-\gamma} \sum_{n=g_\vep(x)}^\infty \ex^{n\,\log \big(\inf_{k\ge 1} (\E [\|\Pi_k\|^\gamma])^{1/k}\big)}\\
&=& c\, \ex^{\log x\big((\alpha-\gamma) + (g_\vep(x)/\log x)(\log h(\gamma)-\log h(\alpha))\big)}\\
&\le & x^{-\delta}\to 0\,,\qquad \xto\,,
\eeao
for some $\delta>0$, depending on $\gamma$ and $\vep$.
\end{proof}
It follows immediately that
\beao
x^\alpha\sum_{n=g_\vep(x)}^\infty \P( \bfv^\top \Pi_n^\top>x)= o(\log x)\,,\qquad \xto\,.
\eeao
\par
As in the univariate case we proceed with an exponential  change of \ms . 
However, the change of measure cannot be done on the marginal distribution but on the transition kernel of the 
Markov chain $({\Pi_n'}^\top\bfu)$. It is indeed a homogeneous  Markov chain as its kernel is given by the expression
\beao
\P(d\bfy| {\Pi_n'}^\top\bfu=\bfx)=\P(\bfX_{n+1}^\top \bfx\in d\bfy)=\P(\bfX^\top \bfx\in d\bfy),
\eeao
that does not depend on $n\ge 1$. The change of kernel is given by
\beao
\P^\alpha(d\bfy\mid {\Pi_n'}^\top\bfu=\bfx)=\ex^{\alpha\, \log (|\bfy|/|\bfx|)}\,\P(\bfX^\top \bfx\in d\bfy).
\eeao
Since this change is difficult to justify we refer the reader for details 
to Buraczewski et al. \cite{buraczewski:damek:guivarch:mentemeier:2014}; see also Section~4.4.6 in the book Buraczewski et al.
\cite{buraczewski:damek:mikosch:2016}. 
In contrast to Kesten's condition in Section~\ref{subsec:kesten},
\cite{buraczewski:damek:guivarch:mentemeier:2014} require a moment condition slightly stronger than in Kesten \cite{kesten:1973}; 
see the discussion on p.~181 in \cite{buraczewski:damek:mikosch:2016}.
This condition is satisfied because we require  
$\E [\|\bfX\|^\alpha (\log \|\bfX\|)^{1+\delta}]<\infty$ for some positive $\delta$.  
\par
Next we notice that $\log |{\Pi_n'}^\top\bfu|$ has the  structure of a Markov random walk associated with $({\Pi_n'}^\top\bfu)$ 
under the change of measure. Indeed, since $|\bfu|=1$ we have the following identity for $n=2$,
\beao\lefteqn{
\P^\alpha(\log |{\Pi_2}'^\top\bfu|\in dx)}\\&=&
\int_{\R^d} \P^\alpha\big(\log |{\Pi_2'}^\top\bfu|\in  d x \mid {\Pi_{1}'}^\top\bfu=\bfy\big)\,\P^\alpha\big( {\Pi_{1}'}^\top\bfu\in d\bfy\big)\\
&=&\int_{\R^d} \ex^{\alpha\, (x-\log|\bfy|)}\,\P \big(\log |{\Pi_2'}^\top\bfu|\in dx\mid {\Pi_{1}'}^\top\bfu=\bfy\big)\,
\ex^{\alpha \,\log|\bfy|}\,\P\big( {\Pi_{1}'}^\top\bfu\in d\bfy\big)\\
&=& \ex^{\alpha\, x}\,\P\big(\log |{\Pi_2'}^\top\bfu|\in dx \big)\,.
\eeao
Thus, the structure of a Markov random walk is recovered thanks to a recursive argument.
\par
From the 
discussion on p.~181 in \cite{buraczewski:damek:mikosch:2016} we also learn that under the aforementioned additional
moment condition,  $(\log \|\Pi_n'\|)$ has a positive drift under $\P^\alpha$, 
i.e., $\log \|{\Pi_n'}\|/n\to h'(\alpha)$ $\P^\alpha$-\as~ On the other hand, by the multiplicative ergodic theorem 
this means that $h'(\alpha)$ is the top Lyapunov exponent under $\P^\alpha$, i.e., $n^{-1} \E^\alpha[\log \|\Pi_n'\|]\to h'(\alpha)$.
\par
As in the univariate case, we have
\beao
\nu(x)&=&x^\alpha \sum_{n=1}^{g_0(x)}\P(|\Pi_n^\top \bfu|>x)=
\int_{0}^\infty \ex^{-\alpha s} d\nu_\alpha(s+\log x)\,,
\eeao
where 
\beao
\nu_\alpha(y)= \sum_{n=1}^{g_0(x)}\P^\alpha(|\Pi_n^\top\bfu|\le y)\,.
\eeao
Since we assume $\E [\|\bfX\|^\alpha (\log \|\bfX\|)^2]<\infty$ we also have  $\E^\alpha [(\log \|\bfX\|)^2]<\infty$. 
By virtue of Theorem~3 in  Hennion \cite{hennion:1997}
we have for any unit vectors $\bfu,\bfv$, 
\beam\label{eq:clthennion}
\sup_z\big|\P^\alpha\big((\log \bfu^\top\Pi_n\bfv-n\,m(\alpha))/\sqrt{n}\le z \big)- \Phi_{0,\sigma^2(\alpha)}(z)\big|\to 0\,,\qquad \nto\,,\nonumber\\ 
\eeam
where $\Phi_{0,\sigma^2(\alpha)}$ is the \df\ of an $N(0,\sigma^2(\alpha))$-distributed \rv , where $\sigma^2(\alpha)\ge 0$ is the limiting variance
which is independent of $\bfu,\bfv$.
An application of Lemma~5.1 in \cite{hennion:1997}
ensures that
\beao
\sup_z\big|\P^\alpha\big((\log |\Pi_n^\top\bfu|-n\,m(\alpha))/\sqrt{n}\le z \big)- \Phi_{0,\sigma^2(\alpha)}(z)\big|\to 0\,,\qquad \nto\,, 
\eeao
For these two results of Hennion one needs the condition \eqref{eq:hennion}.
\par
Now the same arguments as in the univariate case apply to show that 
\beao
\ell(x)&=&\nu_\alpha(\log x)\\& =&\sum_{n=1}^{g_0(x)}\P^\alpha\big((\log |\Pi_n^\top\bfu|-n m(\alpha))/\sqrt{n}\le (\log x-n\,m(\alpha))/\sqrt{n}\big)\\
&=&g_0(x)-\sum_{n=1}^{g_0(x)}\ov \Phi_{0,\sigma^2(\alpha)}\big((\log x-n\,m(\alpha))/\sqrt{n}\big)+ o(\log x)\sim g_0(x)\,. 
\eeao
\par
Now we consider $\P(\bfv^\top \bfY(\bfu)>x)$ for a given $\bfv\in\bbs^{d-1}_+$, $\bfv\ge \bf0$ and $\bfu>0$.
We will mimic the preceding arguments. Since $\bfX \ge 0$ a.s. with non-zero rows ${\Pi_n'}^\top\bfu>0$ a.s. for any $n\ge 1$. Thus, considering the Markov chain $({\Pi_n'}^\top\bfu)_{n\ge 0}$ on the restricted state space $(0,\infty)^d$, one can change the measure according to
\beao
\P^{\alpha,\bfv}(d\bfy\mid {\Pi_n'}^\top\bfu=\bfx)=\ex^{\alpha\, \log (\bfv^\top\bfy/\bfv^\top\bfx)}\,\P(\bfX^\top \bfx\in d\bfy)
\eeao
for any $\bfx$ and $\bfy> \bf0$. Under this change of \ms ,
$\log \bfv^\top{\Pi_n'}^\top\bfu$ has the  structure of a Markov random walk associated with $({\Pi_n'}^\top\bfu)$: 
\beao\lefteqn{
\P^\alpha(\log \bfv^\top{\Pi_2}'^\top\bfu\in dx)}\\&=&
\int_{\R^d} \P^\alpha\big(\log \bfv^\top{\Pi_2'}^\top\bfu\in dx\mid {\Pi_{1}'}^\top\bfu=\bfy\big)\,\P^\alpha\big( {\Pi_{1}'}^\top\bfu\in d\bfy\big)\\
&=&\int_{\R^d} \ex^{\alpha\, (x-\log \bfv^\top\bfy)}\,\P \big(\log \bfv^\top{\Pi_2'}^\top\bfu\in dx\mid {\Pi_{1}'}^\top\bfu=\bfy\big)\,\\
&& \times\ex^{\alpha \,\log (\bfv^\top\bfy/\bfv^\top\bfu)}\,\P\big( {\Pi_{1}'}^\top\bfu\in d\bfy\big)\\
&=& \ex^{\alpha\, (x-\log(\bfv^\top\bfu))}\,\P\big(\log \bfv^\top{\Pi_2'}^\top\bfu\in dx \big)\,.
\eeao
Using the relation  $c|\bfx|\le \bfv^\top\bfx \le c'|\bfx|$ for some $c,c'>0$ and any $\bfx\in (0,\infty)^d$, we have 
\beao
\lim_{n\to \infty} \big(\E [(\bfv^\top\Pi_n'^\top\bfu)^s]\big)^{1/n}=\lim_{n\to \infty} \big(\E [|\Pi_n'^\top\bfu|^s]\big)^{1/n}=\lim_{n\to \infty} \big(\E [\|\Pi_n'\|^s]\big)^{1/n}=h(s).
\eeao
Here, the second identity is stated in Theorem 6.1 of \cite{buraczewski:damek:guivarch:mentemeier:2014}.
From these 3 identities one concludes that the mean of $\bfv^\top\Pi_n'^\top\bfu$  is $h(\alpha)$, the same than the one of $|\Pi_n'^\top\bfu|$ under their respective changes of measure.  For $g_0$ defined as above (independently of $\bfv$) we obtain 
 \beao 
\nu_{\alpha,\bfv}(\log x)=\sum_{n=1}^{g_0(x)} \P^{\alpha,\bfv}(\bfv^\top\Pi_n\bfu\le x)\sim g_0(x).
\eeao
Similar arguments as above and as in the proof of Corollary~\ref{cor:1}  
finish the proof of Theorem~\ref{thm:3} because
\beao
\nu_{\bfv}(x)&:=&\frac{x^\alpha}{(\bfv^\top\bfu)^\alpha} \sum_{n=1}^{g_0(x)}\P(\bfv^\top\Pi_n^\top \bfu>x)=
\int_{0}^\infty \ex^{-\alpha s} d\nu_{\alpha,\bfv}(s+\log x)\,.
\eeao

\appendix
\section{}\label{sec:appA}\setcounter{equation}{0}
In this section we provide an auxiliary result which is the analog of \eqref{eq:3} in the case of positive $X$.
\bpr Assume $\P(X<0)>0$.
If the conditions of Theorem~\ref{thm:kesten} hold then for some $c>0$, 
\beao
\lim_{\xto}x^{\alpha}\,\P\big(\max_{n\ge 1}\Pi'_n> x\big)=c\,.
\eeao
If the conditions of Theorem~\ref{thm:grinc} hold then for some $0<c<c'<\infty$ and for $x$ large enough,
\beao
x^{\alpha}\,\P\big(\max_{n\ge 1}\Pi'_n> x\big)\in[c,c'].
\eeao

\epr
\begin{proof}
Here and in what follows, we interpret $\Pi_0'=1$. We can also work given that $X_n'\neq 0$ as otherwise $\Pi'_t=0$ for all $t\ge n$ and $\max_{1\le t\le n-1}\Pi_t$ has a negligible tail for any fixed $n$.\\

The proof follows the arguments of Goldie \cite{goldie:1991}. One first observes that 
\beam\label{eq:gold1}
\P\big(\max_{n\ge 0}\Pi'_n> x\big)=\E\big[\P\big(\max_{n\ge 0}\Pi'_n> x\mid \bfZ\big)\big]\,,
\eeam
where ${\bf Z}=(Z_n)_{n\ge 1}$ is a Markov chain on $\{-1,1\}$  given by $Z_n=\Pi_n'/|\Pi_n'|$. We define 
\beao
N_0=0\,,\quad N_i=\inf\{k>N_{i-1}: Z_k=1\}\,,\qquad i\ge 1\,\,,\quad I=\{N_i\,,i\ge 1\}\,. 
\eeao
We have
\beam\label{gold2}
\E\big[\P\big(\max_{n\ge 0}\Pi'_n> x\mid {\bf Z}\big)\big]&=&\E\big[\P\big(\max_{n\in I}\Pi'_n> x\mid {\bf Z}\big)\big]\nonumber\\&=&\E\big[\P\big(\max_{n\in I}|\Pi'_n|> x\mid {\bf Z}\big)\big]\nonumber\\
&=&\P(\tau<\infty )\,,
\eeam
where $\tau$ is defined as the stopping time 
\beao
\tau=\inf\Big\{k\in I: \sum_{i=1}^k\log |X_i|>\log x \Big\}\,.
\eeao
We change the \ms\ as follows
\beao
d\P^\alpha(\log |X_n|\le y ,Z_n=\pm 1)= \ex^{\alpha\,y}\,d\P(\log |X_n|\le y, Z_n=\pm 1)\,.
\eeao
We notice that the tilted distribution conditionally on $ {\bf Z}$ is given by 
\beao
d\P^\alpha(\log |X_n|\le y\mid {\bf Z})= \1(Z_n=Z_{n-1})\,\eta_+(dy)+\1(Z_n\neq Z_{n-1})\,\eta_-(dy)\,,
\eeao
where
\beao
\eta_\pm(dy)=\frac{\P^{\alpha}(\pm X>0,\log |X|\in dy )}{\P^\alpha(\pm X>0)}\,.
\eeao
Then we have the identity
\beam\label{eq:kest5}
\P(\tau<\infty )=   \E^{\alpha}\big[ \ex^{-\alpha \sum_{i=1}^\tau\log |X_i|}\1(\tau<\infty) \big].
\eeam
The proof of \eqref{eq:kest5} relies on the fact that $(\ex^{\alpha \sum_{i=1}^n\log |X_i|})$ is a martingale  under $\P$. 
First notice that it is a product of the independent \rv s $|X_i|^\alpha$ given ${\bf Z}$. 
Moreover, noticing that 
\beao
\P(Z_{n}= Z_{n-1})=\P^\alpha(X>0)=\E[|X|^\alpha\1(X>0)]\,,
\eeao 
we have 
\beao
&&\E[\ex^{\alpha  \log |X_n|}\mid {\bf Z}]=\\&&\frac{\E[|X|^\alpha\1(X>0)]}{\P^\alpha(X>0)}\1(Z_{n}= Z_{n-1})+\frac{\E[|X|^\alpha\1(X<0)]}{\P^\alpha(X<0)}\1(Z_{n}\neq Z_{n-1})=1.
\eeao
Then $(\ex^{\alpha \sum_{i=1}^n \log |X_i|})$ constitutes a martingale. Write $\bbf_\tau$ for the $\sigma$-field generated by $\tau$.  Then we have for $n\ge 1$,
\beao\lefteqn{\E^{\alpha}\big[\ex^{-\alpha \sum_{i=1}^\tau\log |X_i|}\1(\tau\le n)\big]}\\&=&\E \big[\ex^{\alpha \sum_{i=1}^n\log |X_i|}\ex^{-\alpha \sum_{i=1}^\tau\log |X_i|}\1(\tau\le n)]\big]\\
&=&\E \big[\E \big[\ex^{\alpha \sum_{i=1}^n\log |X_i|}\ex^{-\alpha \sum_{i=1}^\tau\log |X_i|}\1(\tau\le n)\mid {\bf Z}\big]\big]\\
&=&\E \big[\E \big[\E[\ex^{\alpha \sum_{i=1}^n\log |X_i|}\mid \mathcal F_\tau, {\bf Z}]\,\ex^{-\alpha \sum_{i=1}^\tau\log |X_i|}\1(\tau\le n)\mid {\bf Z}\big]\big]\\
&=&\E\big[\E [ \ex^{\alpha \sum_{i=1}^\tau\log |X_i|} \ex^{-\alpha \sum_{i=1}^\tau\log |X_i|}\1(\tau\le n)\mid {\bf Z}]\big]\\
&=&\P(\tau\le n ).
\eeao
Now \eqref{eq:kest5} follows by letting $n\to \infty$.
\par
In view of \eqref{eq:gold1}--\eqref{eq:kest5} it suffices to study the \asy\ behavior of $x^\alpha\,\P(\tau<\infty)$. 
From \eqref{eq:kest5} we conclude that
\beao
x^{\alpha}\,\P(\tau<\infty )&=&   \E^{\alpha}[\E^{\alpha}[\ex^{-\alpha(\sum_{i=1}^\tau\log |X_i|-\log x)}\1(\tau<\infty)\mid {\bf Z}]].
\eeao
We have the \rep\ 
\beao
\sum_{i=1}^\tau\log |X_i|=\sum_{i=1}^{\tilde\tau} \sum_{j=N_{i-1}+1}^{N_i} \log |X_i| =:\sum_{i=1}^{\tilde \tau}W_i
\eeao
where $(W_i)$ is iid and has the common tilted \ds\ $\eta$ given explicitly in (9.11) of Goldie \cite{goldie:1991}. Here   
\beao
\tilde \tau=\inf \{t\ge 1:\sum_{i=1}^{t}W_i \ge \log x \}\,.
\eeao 
It is shown in equations (9.12) and (9.13) of \cite{goldie:1991} that $  \E^\alpha[W_1]=2\,m(\alpha)>0$.\\
\par
Assume now that the conditions of Theorem \ref{thm:kesten} are satisfied and that the distribution of $\log|X|$ given $X\neq 0$ is non-arithmetic. Then $\eta$ inherits non-arithmeticity from $\log |X|$ and the overshoot $B(\infty)$ of the random walk 
associated  with $\eta$ is well defined (positive drift) and
\begin{align*}
x^{\alpha}\,\P(\tau<\infty )&=   \E^{\alpha}\big[\E^{\alpha}[\ex^{-\alpha(\sum_{i=1}^{\tilde \tau}W_i-\log x)}\1(\tilde \tau<\infty)\mid {\bf Z}]\big]\\
&=\E^{\alpha}[\ex^{-\alpha(\sum_{i=1}^{\tilde \tau}W_i-\log x)} ]\\
&\to \E^{\alpha}[\ex^{-\alpha B(\infty)} ]=\frac{1-h(0,\infty)}{2\,\alpha \,m(\alpha)}>0\,.
\end{align*}
Here $H$ is the ladder height (defective) distribution of the random walk associated with the distribution of $W_1$ under $\P$ (negative drift) and the positivity of the constant on the \rhs\ follows by calculations similar to 
Theorem 5.3, Section XIII, in  Asmussen \cite{asmussen:2003}. 
This finishes the proof of the first assertion.\\

Considering  that the assumptions of Theorem \ref{thm:grinc} are satisfied then the distribution of the overshoot $B(\infty)$ of the random walk 
associated  with $\eta$ is well defined only on some lattice span only. Then $x^{\alpha}\,\P(\tau<\infty )$ converges when $\xto$ on this lattice span only (with a different limiting positive constant), see Remark 5.3, Section XIII, in  Asmussen \cite{asmussen:2003}. The second assertion follows by using the monotonicity of the tail probability.
\end{proof}
\bre
The distribution of $W_1=\log \tilde X_1$ under $\P$ given in \eqref{eq:distr} heavily depends on the quantities $p$ and $q$. 
Therefore we are not convinced that 1. $H$ ever coincides with the ladder height \ds\ of the random walk with step sizes
$(\log |X_n|)$ and 2. the relation
\beao
2\, x^{\alpha}\,\P\big(\max_{k\ge 0}\Pi'_k> x\big)\sim x^{\alpha}\,\P\big(\max_{k\ge 0}|\Pi'_k|> x\big),\qquad x\to\infty,
\eeao
holds. It is clearly not the case when $p=0$ ($N_1=2$) as shown by an inspection of Spitzer's formula; see \cite{feller:1971}, p. 416.
\ere
\subsection*{Acknowledgments} { This research was started in January 2017 when Moh\-sen Rezapour visited the Department
of Mathematics at the University of Copenhagen. He takes pleasure to thank his Danish colleages for their hospitality.
The fact that $Y$ is heavy-tailed was communicated to us by Paul Embrechts who reacted to a corresponding question
by Professor Didier Sornette. From Professor Sornette we learned recently that this topic was actually developed in the research 
team of Professor Misako Takayasu and Professor Hideki Takayasu with their PhD student Arthur Matsuo at Tokyo Tech.
We are grateful to Professor Sornette for posing this problem and for discussions on the topic.}

\end{document}